\documentclass[12pt]{article}
\usepackage{mathrsfs}
\usepackage{amssymb,amsfonts}
\def\<{\langle}
\def\>{\rangle}
\def\a{\alpha}

\def\bu{\bullet}
\def\bm{\boxminus}
\def\bt{\boxtimes}

\def\ci{\circ}

\def\di{\diamondsuit}
\def\btr{\blacktriangleright}
\def\btl{\blacktriangleleft}
\def\D{\Delta}

\def\G{\Gamma}

\def\i{\iota}

\def\lr{\longrightarrow}

\def\o{\otimes}

\def\U{\Upsilon}

\def\t{\tau}

\def\si{\sigma}
\def\v{\varepsilon}
\def\wh{\widehat}

\date{}
\begin{document}
\renewcommand{\baselinestretch}{1.2}
\renewcommand{\arraystretch}{1.0}
\title{\bf   Twisted  Algebras of Multiplier Hopf ($^*$-)algebras}
\date{}
\author {{\bf  Shuanhong Wang }\\
{\small Department of Mathematics, Southeast University}\\
{\small Nanjing, Jiangsu 210096, P. R. of China}\\
{\small E-mail:   shuanhwang@seu.edu.cn}}
 \maketitle
\begin{center}
\begin{minipage}{12.cm}
\begin{center}{\bf ABSTRACT}\end{center}

In this paper we study twisted  algebras of multiplier
 Hopf ($^*$-)algebras  which generalize all kinds of
 smash products such as generalized smash products,
 twisted smash products, diagonal crossed products,
 L-R-smash products, two-sided crossed products and
  two-sided smash products for the ordinary
  Hopf algebras appeared in [P-O].



\end{minipage}
\end{center}

\section*{1. Introduction }
\def\theequation{0. \arabic{equation}}
\setcounter{equation} {0} \hskip\parindent

In [P-O], Panaite  and  Van Oystaeyen  introduced a more general version of the so-called L-R-smash product
 for (quasi-)Hopf algebras and studied its relations with
 other kinds of crossed products (two-sided smash and crossed product and diagonal crossed product),
  In order to find the method of constructing multiplier
 Hopf ($^*$-)algebras, in this paper  we will try to generalize these results to the situation for multiplier Hopf algebras introduced
  in [VD1-VD3].
 \\

This paper is organized as follows. In Section 2 of this paper  we will study the notions of a generalized
 twisted tensor product, a generalized twisted smash product
  and a generalized L-R-smash product. We also study some
  isomorphisms between them.
 In Section 3 we focus on two-times twisted tensor products for multiplier Hopf algebras.
 Finally, in Section 4 we will introduce the notion of a Long module algebra and study the
 the twisted product of the given multiplier Hopf algebras.
\\

\section*{1. Preliminaries }
\def\theequation{1. \arabic{equation}}
\setcounter{equation} {0} \hskip\parindent

Let $A$ be  an associative algebra with a nondegenerate product.
 If $A$ has a unit element, this requirement is automatically satisfied.
 In the general case, we only suppose the existence of local units in the following sense.
  Let $\{e_1, \cdots  , e_n\}$ be a finite set of elements in $A$.
   Then there exist elements $e, f \in A$ so that
   $ee_i = e_i = e_i f$ for $i = 1, \cdots , n$.
\\

 A multiplier $m=(m_1, m_2)$ of the algebra $A$ is a pair of linear mappings in $End_k(A)$
  such that $m_2(a)b = a m_1(b)$ for all $a, b \in A$. The set of multipliers
  of $A$ is denoted by $M(A)$. It is a unital algebra which contains $A$
   as essential ideal through the embedding $a\hookrightarrow  (a \cdot, \, \cdot a)$.
   Therefore $m\cdot a = (m_1(a)\cdot , \, \cdot m_1(a))\equiv m_1(a)$
    and $a \cdot m = (m_2(a)\cdot ,\,  \cdot m_2(a))\equiv m_2(a)$ for all
    $m\in  M(A)$ and $a\in  A$. Hence we will frequently use the identification
    $a\cdot m = m_2(a)$ and $m\cdot a = m_1(a)$. If $A$ is unital
   then $A = M(A)$. If $A$ is a $*$-algebra then $M(A)$ is a  $*$-algebra through
   $m^* = (m_2^* ,\,  m_1^* )$ where $\psi ^*(a) := \psi (a^*)^*$ for any
   $a\in A$, $\psi \in End_k(A)$. Since the multiplication of $A$
    is supposed to be nondegenerate a multiplier
    $m=(m_1, m_2)$ of $A$ is uniquely determined by its first or second component.
 For a tensor product of two algebras $A$ and $B$ one obtains the canonical algebra embeddings
  $A\o B \hookrightarrow  M(A)\o M(B)\hookrightarrow M(A\o B)$.
\\

Let $A$ be an algebra with a non-degenerate product. An algebra
 homomorphism $\D : A\lr M(A\o A)$ is called a comultiplication
 on $A$ if

 (i) $T_1(x\o y)=\D (x)(1\o y)\in A\o A$ and  $T_2(x\o y)=(x\o 1)\D (y)\in A\o A$
 for all $x, y\in A$ (some times we call this condition a "coving theory");

 (ii) $(T_2\o id)(id\o T_1)=(id\o T_1)(T_2\o id)$ on $A\o A\o A$.
 \\

 We call $A$ with a comultiplication $\D $ a multiplier Hopf algebra if the linear
 maps $T_1, T_2$ on $A\o A$ are bijective. Moreover, we call $A$ regular if $\t \D $,
 where $\t $ is the flip, is again a comultiplication such that $(A, \t \D)$ is also a
 multiplier Hopf algebra. If $A$ is a $*$-algebra, we call $\D $ a comultiplication
 if it is also a $*$-homomorphism. A multiplier Hopf $*$-algebra is a $*$-algebra
  with a comultiplication, making it into a multiplier Hopf algebra.
 \\

An algebraic quantum group is a regular multiplier Hopf algebra with non-trivial invariant functionals
 (i.e., integrals).
\\

Notice that  we can use the Sweedler's notation for $\D (x)$ when $x\in A$. The problem is
 that $\D (x)$ is not in $A\o A$ in general. By covering theory, we know however that $\D (x)(1 \o y) \in A\o  A$
 for all $x, y\in A$. We can write (cf.  [Dr-VD-Z, VD3]) $\D (x)=\sum x_1\o x_2$ for this.
 Now, we know that there is an element $e \in A$ such that $y = ey$ and we can think
 of $\D (x)=\sum x_1\o x_2$ to stand for $\D (x)(1 \o e)$. Of course, this is still dependent on $y$.
 But we know that for several elements $y$, we can use the same $e$.
\\

 Let $A$ be a regular multiplier Hopf algebra
 (i.e. with a bijective antipode). By a left $A$-module $V$,
  we always mean a unital module (sometimes, we also say that
  $\rhd$ is a left action of $A$ on $V$). This means that $A\rhd  V= V$.
   For all $v\in V$, we have an element $e\in A$ such that
   $e\rhd v =v$. Observe that for a unital algebra $A$,
   this condition is automatically satisfied. Similarly, a right $A$-module
    can be defined.
\\

Let $Q$ denote a regular multiplier Hopf algebra and $A$ an associative algebra with or without
 identity.  Assume that $\rhd$ is a left action of $Q$ on $A$
  making $A$ into a left $Q$-module algebra, i.e., we have
  $x\rhd (ab)=\sum (x_1\rhd a)(x_2\rhd b)$
  for all $x\in Q$ and $a, b\in A$. In this formula, $x_1$ is covered by $a$ (through the action)
  and  $x_2$ is covered by $b$ (through the action).
   Similarly,  assume that $\lhd$ is a right action of $Q$ on $A$
  making $A$ into a right $Q$-module algebra, i.e., one has
  $(ab) \lhd x =\sum ( a\lhd x_1)( b\lhd x_2)$,
  for all $x\in Q$ and $a, b\in A$.
   Then $A$ is called a $Q$-bimodule algebra, if it
 is a unital left $Q$-module and a unital right $Q$-module
 such that $(x\rhd a)\lhd y=x\rhd (a\lhd y)$,
  for all $x, y\in Q$ and $a\in A$.
\\

 Suppose that $\G$ is a left coaction of $Q$ on $A$ making
  $A$ a left $Q$-comodule algebra.  More precisely, $\G : A\lr M(Q\o A)$
  is an injective homomorphism so that for all
  $x\in Q$ and $a\in A$, we have
  $\G (a)(x\o 1)$ and $(x\o 1)\G (a)$ are in $Q\o A$.
  We use the Sweedler notation for these expressions, e.g.
  $\G (a)(x\o 1) = \sum a_{(-1)}x\o a_0$.
   We say that the multiplier $\G (a)$ is covered by
   $x\o 1$. Now the further requirement makes sense:
    $(\i \o \G )\G(a)=(\D \o \i)\G (a)$ for all $a\in A$.
    In fact, the expression  $(\i \o \G )\G(a)=\sum a_{(-2)}\o a_{(-1)}\o a_0$
    makes sense when it is understood that $\D (a_{(-1)})$  is replaced by $a_{(-2)}\o a_{(-1)}$. We have more or less
  the same rules for covering. We need to cover the factor $ a_{(-1)}$ and possibly also $a_{(-2)}$ (and so on) by
  elements in $A$, left or right. In this case however, one cannot cover the factor $ a_0$.
     Similarly,  the right coaction  $\U$ of $L$ on $A$ can be defined.
    Then $A$ is called a $Q$-$L$-bicomudule algebra
    if
 $$
 (1\o 1\o  x)((\i \o \U)(\G (a)(y\o 1)))=((\G \o \i )((1\o x)\U (a)))(y\o 1)
 $$
for all $x, y\in Q$ and $a\in A$.
\\

Let $Q, L$ be multiplier Hopf algebras. Denote the category of all
  left $Q$-module algebras by ${}_Q{\cal MA}$, denote by the category
  of all right $L$-module algebras ${\cal MA}_L$, and
  denote by ${}_Q{\cal MA}_L$ the category of all $Q$-$L$-bimodule algebras.
  Similarly,  denote the category of all
  left $Q$-comodule algebras by ${}^Q{\cal MA}$, denote by the category
  of all right $L$-module algebras ${\cal MA}^L$, and
  denote by ${}^Q{\cal MA}^L$ the category of all $Q$-$L$-bicomodule algebras.
\\

Let $Q$ and $\wh{Q}$ be a dual pair of algebraic quantum groups over $\mathbb{C}$. By Pontryagin duality,
 identifying $\wh{\wh{Q}}=Q$, we have a natural
 dual pairing $Q\o \wh{Q}\lr \mathbb{C}$ is written as
 $$
\< x, \wh {x}\>\in \mathbb{C}, \quad x\in Q, \, \, {\mbox{and}}\, \, \wh{x}\in \wh{Q}.
 $$
Note that the antipode of $Q^{op}$ and $Q^{cop}$ is given by
 $S^{-1}$ and the antipode of $Q^{op,cop}$ by $S$. Also,
 $\wh{Q^{op}}=\wh{Q}^{cop}$,  $\wh{Q^{cop}}=\wh{Q}^{op}$
 and  $\wh{Q^{op, cop}}=\wh{Q}^{cop,op}$.
\\

Obviously, after a permutation of tensor factors $A\o Q\longleftrightarrow Q\o A$
 a left coaction of $Q$ may always be viewed as a right coaction by $Q^{cop}$ and vice versa.
\\

Next, we recall that there is a one-to-one correspondence between right (left) coactions of $Q$
 on $V$ and left (right) actions, respectively, of $\wh{Q}$ on $V$
  given for $\wh{x}\in \wh{Q}$ and $v\in V$ by
  $$ \wh{x} \rhd v=(\i \o \wh{x})\U (v) $$ and
  $$ v \lhd \wh{x}=(\wh{x}\o \i ) \G (v).$$

 As a particular example we recall the case $V=Q$ with
  $\G =\U =\D $. In this case we denote the associated left and right actions
  of $\wh{x}\in \wh{Q}$ on $y\in Q$ by $\wh{x}\btr y$ and $y\btl \wh{x}$, respectively,
  see [De].
  \\

We refer to [VD1-VD3]  for the theory of multiplier Hopf algebras.
 For the use of the Sweedler notation in this setting, we refer to [Dr-VD] and [Dr-VD-Z].
  For pairings of multiplier Hopf algebras, the main
 reference is [Dr-VD].  An important reference for this paper is of course [P-O] where
 the theory of a more general version of the so-called L-R-smash product
 for (quasi-)Hopf algebras  is developed. Finally, for some
 related constructions with multiplier Hopf algebras, see [De, De-VD-W, De-VD] and [VD-VK, W, W-L].

\section*{2. Twisted tensor products  }

In this section, we will study the notions of a generalized
 twisted tensor product, a generalized twisted smash product
  and a generalized L-R-smash product. We also study some
  isomorphisms between them.
\\

{\sl  2.1 Generalized twisted tensor products}
\\
\def\theequation{2. 1. \arabic{equation}}
\setcounter{equation}{0} \hskip\parindent

Let $A, B$ be algebras  and suppose that there are two
 $k$-linear maps $R: B\o A\lr A\o B$ and $T: A\o B\lr A\o B$
 such that for $a, c\in A, b, d\in B$
\begin{eqnarray*}
(2.1) && a_R\o (bd)_R=a_{Rr}\o b_{r}d_R \, \Leftrightarrow \, R(m_B\o id _A)=(id_A\o m_B)R_{12}R_{23}; \\
(2.2) && (ac)_R\o b_R=a_{R}c_{r} \o b_{Rr} \, \Leftrightarrow \, R(id_B\o m_A)=(m_A\o id_B)R_{23}R_{12};\\
(2.3) &&  a_T\o (bd)_T=a_{Tt}\o b_{T}d_t  \, \Leftrightarrow \, T(id_A\o m_B)=(id_A\o m_B)T_{13}T_{12}; \\
(2.4) && (ac)_T\o b_T=a_{T}c_{t} \o b_{tT}  \, \Leftrightarrow \, T(m_A\o id_B)=(m_A\o id_B)T_{13}T_{23};\\
(2.5) &&  a_T\o c_R \o d_{TR}= a_T\o c_R \o d_{RT}\, \Leftrightarrow \, R_{23}T_{12}=T_{13}R_{23}.
\end{eqnarray*}
Here we write  $R(b\o a)=a_R\o b_R=a_r\o b_r$ and
 $T(a\o b)=a_T\o b_T=a_t\o b_t$ for all $a\in A$ and $b\in B$,
  and  $R_{12}=(R\o \i)\in Hom _k(B\o A\o X, A\o B\o X)$,
   $R_{23}=(\i\o R)\in Hom _k(X\o B\o A, X\o A\o B)$,
      $T_{12}=(T\o \i)\in Hom _k(A\o B\o X, A\o B\o X)$,
  $T_{13}\in Hom _k(A\o X\o B, A\o X\o B)$ works on
   the first and the third componets as $T$ and
  $T_{23}=(\i\o T)\in Hom _k(X\o A\o B, X\o A\o B)$ for any $k$-vector space
  $X$.
\\

One can define a {\sl generalized L-R-twisted tensor product} $A\# _{RT}B=A\o B$
 as spaces with a new multiplication given by the formula:
\begin{equation}
m_{A\# _{RT}B}=(m_A\o m_B)\circ  T_{14}\circ R_{23}\quad  \mbox
{or}\quad (a\# b)(c\# d)=a_Tc_R\# b_Rd_T
\end{equation}
for all $a, c\in A$ and $b, d\in B$.  Here,  $T_{14}\in End _k(A\o B\o A\o B)$ works on
   the first and the forth componets as $T$.
\\

{\bf Proposition 2.1.1.} The above multiplication on $A\# _{RT}B=A\o B$
 is associative.

{\bf Proof.}  For all $a, c, x\in A$ and $b, d, y\in B$. We compute
\begin{eqnarray*}
[(a\# b)(c\# d)](x\# y)&=&(a_Tc_R)_tx_r\# (b_Rd_T)_ry_t\\
&\stackrel {(2.4)}{=}&a_{Tt}c_{Rt'}x_r\# (b_Rd_T)_ry_{t't}\\
&\stackrel {(2.1)}{=}&a_{Tt}c_{Rt'}x_{rr'}\# b_{Rr'}d_{Tr}y_{t't},
\end{eqnarray*}
and
\begin{eqnarray*}
(a\# b)[(c\# d)(x\# y)]&=& a_T(c_tx_r)_R\# b_R(d_ry_t)_T\\
&\stackrel {(2.3)}{=}& a_{Tt'}(c_tx_r)_R\# b_R d_{rT}y_{tt'}\\
&\stackrel {(2.2)}{=}& a_{Tt'} c_{tR}x_{rR'}\# b_{RR'} d_{rT}y_{tt'}\\
&\stackrel {(2.5)}{=}& a_{Tt'} c_{Rt}x_{rR'}\# b_{RR'} d_{rT}y_{tt'}.
\end{eqnarray*}

 This completes the proof. $\blacksquare$
 \\

 {\bf Remark 2.1.2.} When $T$ is an identity map, then
 a generalized L-R-twisted tensor product $A\# _{RT}B=A\o B$
 becomes a twisted tensor product $A\# _RB=A\o B$.
 As a space, $A\# _RB=A\o B$ with the following product
  defined by the formula:
$$
m_{A\# _RB}=(m_A\o m_B)\circ  R_{23}\quad  \mbox
{or}\quad (a\# b)(c\# d)=ac_R\o b_Rd
$$
for all $a, c\in A$ and $b, d\in B$.
 Here, $R_{23}=\i_A\o R\o \i_B: A\o B\o A\o B\lr A\o A\o B\o B$.
\\

  If $A$ and and $B$ have a unit $1$, we suppose that
  $R$ satisfies $R(1\o a)=a\o 1$ and $R(b\o 1)=1\o b$ for all $a\in A$ and $b\in B$. Then
  $1\o 1$ is the unit of $A\# _RB$. Furthermore, the maps
  $A\lr A\# _RB: a\mapsto a\o 1$ and  $B\lr A\# _RB: b\mapsto 1\o b$
  are algebra embeddings. For more details on the above results, we refer to [VD-VK].
   \\

{\bf Definition 2.1.3.}  Let $A$ and $B$ be  two algebras without unit,
 but with non-degenerate products. Let  $R: B\o A\lr A\o B$ and $T: A\o B\lr A\o B$
  be two $k$-linear maps. Then one can define the following two
  $k$-linear maps $R*T$ and $T*R$ from $A\o B\o A\o B $  to $ A\o B$:
\begin{eqnarray*}
&&(R*T)(a\o x\o b\o y)=(a b_R)_T\o (x y_T)_R;\\
&&(T*R)(b\o y\o a\o x)=(b_Ta)_R\o (y_Rx)_T
\end{eqnarray*}
for all $a, b\in A$  and $ x, y\in B$.
\\

{\bf Proposition 2.1.4.}  Let $A$ and $B$ be two algebras without unit,
 but with non-degenerate products. If $R$ and $T$ are bijective
 such that (1): $(T*R)(x\o y\o a\o b)=0$ implies that $x\o y=0$
 and (2): $(R*T)(a\o b\o x\o y)=0$ implies that $x\o y=0$
 for all $a, x\in A$ and $b, y\in B$,
  then the product in $A\# _{RT}B$ is non-degenerate.

 {\bf Proof.} Suppose that there is an element $\sum a_i\# b_i\in A\# _{RT}B$
 such that $(\sum a_i\# b_i)$ $(a\# b)=0$ for all $a\in A, b\in B$. Then we have that
 $$
 \sum a_{iT} a_r\# b_{ir} b_T=0,
 $$
 this implies (by Eq.(2.2)) that
 $$
 \sum (a_{iTR^{-1}} a)_R\# b_{iR^{-1}R} b_T=0,
 $$
 and by Eq.(2.3) and Eq.(2.5), one has that
 $$
 \sum (a_{it^{-1}R^{-1}T} a)_R\# (b_{it^{-1}R^{-1}R} b)_T=0.
 $$

 By the assumption (1), we have
 $$
 \sum a_{it^{-1}R^{-1}}\# b_{it^{-1}R^{-1}}=0.
 $$

 Thus, we obtain that  $\sum a_i\# b_i=0$.

  Similarly, by Eq.(2.1), (2.4), (2.5) and the assumption (2),
  one can show that
 $(a\# b)(\sum a_i\# b_i)=0$ for all $a\in A, b\in B$ implies that $\sum a_i\# b_i=0$. $\blacksquare$
 \\

{\bf Corollary 2.1.5.}  Let $A$ and $B$ be two algebras without unit,
 but with non-degenerate products.

 (1) ([De, Proposition 1.1]) If $R$ is bijective, then the product in
 $A\# _RB$ is non-degenerate.

(2) If $T$ is bijective, then the product in
 $A\# _TB$ is non-degenerate.

{\bf Proof.}  (1) In this case, $T=id _{A\o B}$. Thus the conditions (1): $(x_Ta)_R\o (y_Rb)_T=0$ implies that $x\o y=0$ and (2): $(a x_R)_T\o (b y_T)_R=0$ implies that $x\o y=0$ for all $a, x\in A$ and $b, y\in B$,
 become that (1): $(xa)_R\o y_Rb=0$ implies that $x\o y=0$
 and (2): $ax_R\o (b y)_R=0$ implies that $x\o y=0$ for all $a, x\in A$ and $b, y\in B$.
 We now explain how these two conditions hold. In fact, for the condition (1), since the product on $B$ is non-degenerate, $(xa)_R\o y_Rb=0$ implies that
 $(xa)_R\o y_R=0$ and thus $ x\o y=0$ since the product on $B$ is non-degenerate. Similarly,
 the condition (2) holds too.

 (2) In this case, $R=id _{A\o B}$. The check is similar to the one in the part (1). $\blacksquare$
 \\

{\bf Theorem 2.1.6.} Let $A$ and $B$ be Hopf algebras with two
 $k$-linear maps $R: B\o A\lr A\o B$ and $T: A\o B\lr A\o B$
 such that Eq.(2.1)-(2.5) hold and
 \begin{eqnarray*}
(2.6) && a_R\o 1_R=a\o 1,\, \,   1_R\o b_R=1\o b; \\
(2.7) && a_T\o 1_T=a\o 1, \, \,  1_T\o b_T=1\o b; \\
(2.8) && (\v _A\o \v_B)R=(\v _B\o \v_A), \, \,   (\v _A\o \v_B)T=(\v _A\o \v_B); \\
(2.9) && (id_A\o \t \o id_B)(\D _A\o \D _B)R=(R\o R)(id_B\o \t \o id_A)(\D _B\o \D _A);\\
(2.10) && (id_A\o \t \o id_B)(\D _A\o \D _B)T=(T\o T)(id_A\o \t \o id_B)(\D _A\o \D _B).
\end{eqnarray*}
Here we write  $R(b\o a)=a_R\o b_R=a_r\o b_r$ and
 $T(a\o b)=a_T\o b_T=a_t\o b_t$ for all $a\in A$ and $b\in B$.

Then  $(A\# _{RT}B, \D _{RT}, \v _{RT}, S_{RT})$ is a bialgebra, where the comultiplication,
 $\D _{RT}$ and $\v _{RT}$ are given as

\begin{eqnarray*}
&&(a\# b)(c\# d)=a_Tc_R\# b_Rd_T, \\
&& \D _{RT}=(id_A\o \t \o id_B)(\D _A\o \D _B),\\
&& \v _{RT}=\v _A\o \v _B,\\
\end{eqnarray*}
for all $a, c\in A$ and $b, d\in B$.

Furthermore, if
\begin{eqnarray*}
(2.11) && \sum (S(a_1)_{T}a_{2})_R\# (S(b_1)_{R}b_{2})_T= \v _{RT}(a\# b)(1_A\# 1_B),\\
(2.12) && \sum (a_1S(a_2)_{R})_T\# (b_1S(b_2)_T)_R= \v _{RT}(a\# b)(1_A\# 1_B),
\end{eqnarray*}
 then   $(A\# _{RT}B, \D _{RT}, \v _{RT}, S_{RT})$ is a Hopf algebra with $S_{RT}$ given by
$$
S_{RT}= T R(S_B\o S_A)\t.
$$

{\bf Proof.} From Proposition 2.1.1, we know that
 $A\# _{RT}B$ is an algebra with unity $1_A\# 1_B$. Clearly
 $\D _{RT}$ is coassociative and $\v _{RT}$ is counity.

 Take $a, c\in A$ and $b, d\in B$. We now prove that
 $ \D _{RT}$ is a homomorphism.
 \begin{eqnarray*}
&&\D _{RT}[(a\# b)(c\# d)]\\
&& =\D _{RT}(a_Tc_R\# b_Rd_T) \\
&& =(id_A\o \t \o id_B)(\D _A\o \D _B)(a_Tc_R\# b_Rd_T) \\
&& =\sum (a_Tc_R)_1\# (b_Rd_T)_1 \o (a_Tc_R)_2\# (b_Rd_T)_2\\
&& =\sum a_{T1}c_{R1}\# b_{R1}d_{T1} \o a_{T2}c_{R2}\# b_{R2}d_{T2}\\
&& \stackrel {(2.9)(2.10)}{=}\sum a_{1T}c_{1R}\# b_{1R}d_{1T} \o a_{2T}c_{2R}\# b_{2R}d_{2T}\\
&& =\D _{RT}(a_Tc_R\# b_Rd_T)\D _{RT}(a_Tc_R\# b_Rd_T).
\end{eqnarray*}

From  the counitary property (2.8) we obtain that $\v _{RT}$ is a homomorphism
 on $A\# _{RT}B$.  Therefore, $(A\# _{RT}B, \D _{RT}, \v _{RT}, S_{RT})$ is a bialgebra.

 We show that $m_{A\# _{RT}B}(S_{RT}\o id_{A\o B})(\D _{RT}(a\# b))=
 \v _{RT}(a\# b)(1_A\# 1_B)$.
 \begin{eqnarray*}
&&m_{A\# _{RT}B}(S_{RT}\o id_{A\o B})(\D _{RT}(a\# b))\\
&&=\sum S(a_1)_{RTt}a_{2r}\# S(b_1)_{RTr}b_{2t}\\
&&\stackrel {(2.5)}{=}\sum S(a_1)_{RTt}a_{2r}\# S(b_1)_{RrT}b_{2t}\\
&&\stackrel {(2.3)}{=}\sum S(a_1)_{RT}a_{2r}\# (S(b_1)_{Rr}b_{2})_T\\
&&\stackrel {(2.5)}{=}\sum S(a_1)_{TR}a_{2r}\# (S(b_1)_{Rr}b_{2})_T\\
&&\stackrel {(2.2)}{=}\sum (S(a_1)_{T}a_{2})_R\# (S(b_1)_{R}b_{2})_T\\
&&\stackrel {(2.11)}{=} \v _{RT}(a\# b)(1_A\# 1_B).
\end{eqnarray*}

Similarly, by Eq.(2.1), (2.4), (2.5) and Eq.(2.12), we can get
\begin{eqnarray*}
&&m_{A\# _{RT}B}( id_{A\o B}\o S_{RT})(\D _{RT}(a\# b))\\
&&=\sum (a_1S(a_2)_{R})_T\# (b_1S(b_2)_T)_R\\
&&= \v _{RT}(a\# b)(1_A\# 1_B).
\end{eqnarray*}

Thus, $(A\# _{RT}B, \D _{RT}, \v _{RT}, S_{RT})$ is a Hopf algebra. $\blacksquare$
\\
\\

{\bf Remark 2.1.7.} In the $(A\# _{RT}B, \D _{RT}, \v _{RT}, S_{RT})$. We have

(1) when $T$  is the identity map, i.e, $T=id_{A\o B}$,
 then the conditions (2.11) and (2.12) hold. In fact,
 we have
 \begin{eqnarray*}
&& \sum (S(a_1)_{T}a_{2})_R\# (S(b_1)_{R}b_{2})_T\\
&&= \sum (S(a_1)a_{2})_R\# S(b_1)_{R}b_{2}\\
&&= \sum \v (a)1_R\# S(b_1)_{R}b_{2}\\
&&\stackrel {(2.6)}{=} \sum \v (a)1\# S(b_1)b_{2}\\
&&= \v _{RT}(a\# b)(1_A\# 1_B),
\end{eqnarray*}
and so Eq.(2.11) is obtained. Similarly,
\begin{eqnarray*}
\sum (a_1S(a_2)_{R})_T\# (b_1S(b_2)_T)_R= \v _{RT}(a\# b)(1_A\# 1_B).
\end{eqnarray*}

In this case, we have already obtained the result of [De, Theorem 2.1].

(2)  when  $R$ is  the flip map, i.e., $\t$,  then the conditions (2.11) and (2.12) hold
 too.

(3) When $R$ and $T$ are not trivial, we will give the example (see Example 2.3.1 (3)).
 \\

 {\bf Proposition 2.1.8.} Let $A$ and $B$ be Hopf $*$-algebras with two
 $k$-linear maps $R: B\o A\lr A\o B$ and $T: A\o B\lr A\o B$
 such that $R$ and $T$ satisfy all the conditions of Theorem 2.0.6. If furthermore
 \begin{eqnarray*}
(2.13) && (R(*_B\o *_A)\t)^2=id _A\o id_B,\\
(2.14) && (T(*_A\o *_B))^2=id _A\o id_B,
\end{eqnarray*}
 then   $(, \D _{RT}, \v _{RT}, S_{RT})$ is a Hopf $*$-algebra

{\bf Proof.}  By Theorem 2.1.6, one has that $(A\# _{RT}B, \D _{RT}, \v _{RT}, S_{RT})$ is a Hopf algebra.
 By the conditions that $R(*_B\o *_A)\t$ and $T(*_A\o *_B)$ are involutions on $A\o B$, we have
 that $A\# _{RT}B$ is a $*$-algebra when $(a\# b)^*=TR(b^*\o a^*)$. In what follows, we check that $\D_{RT}$
  is a $*$-homomorphism on   $A\# B$.
\begin{eqnarray*}
\D _{RT}((a\# b)^*)&=& \D _{RT}(TR(b^*\o a^*))\\
&&=(id_A\o \t \o id_B)(\D _A\o \D _B)TR(b^*\o a^*)\\
&&\stackrel {(2.10)(2.9)}{=}(TR\o TR)(id_B\o \t \o id_A)(\D _B(b^*)\o \D _B(a^*))\\
&&=TR(b^*{}_1\# a^*{}_1)\o TR(b^*{}_2\# a^*{}_2)\\
&&=(a_1\# b_1)^*\o (a_2\# b_2)^*\\
&&=(\D _{RT}(a\# b))^*.
\end{eqnarray*}
This finishes the proof. $\blacksquare$
\\

{\bf Proposition 2.1.9.} Let $m \in  M(A), n\in M(B)$, $M \in  M(A\o A)$ and $N\in M(B\o B)$
 be multipliers. Then

 (1) we have that $m\# n\in M(A\# _{RT}B)$ where $m\# n$
  is defined by
\begin{eqnarray*}
&& (m\# n)_1(a\# x)=m_1(a_{T^{-1}R^{-1}R})_T\# n_1(x_{T^{-1}R^{-1}T})_R; \\
&& (m\# n)_2(a\# x)=m_2(a_{T^{-1}R^{-1}T})_R\# n_2(x_{T^{-1}R^{-1}R})_T
\end{eqnarray*}
for all $a\in A$ and $x\in B$.

(2) One has that  $M\# N\in M((A\# _{RT}B)\o (A\# _{RT}B))$ where $M\# N$
  is defined as follows:
\begin{eqnarray*}
&& (M\# N)_1((a\# x)\o (b\# y))= M^1_1(a_{T^{-1}R^{-1}R})_T\# N^1_1(x_{T^{-1}R^{-1}T})_R \\
&& \quad \quad \quad \quad\quad \quad \quad \quad \quad \quad\quad
\quad \quad \quad \quad\quad  \o \,  M^2_1(b_{t^{-1}r^{-1}r})_t\# N^2_1(y_{t^{-1}r^{-1}t})_r; \\
&& (M\# N)_2((a\# x)\o (b\# y))= M^1_2(a_{T^{-1}R^{-1}T})_R \# N^1_2(x_{T^{-1}R^{-1}R})_T \\
&& \quad \quad \quad \quad\quad \quad \quad \quad \quad \quad\quad
\quad \quad \quad \quad\quad  \o \, M^2_2(b_{t^{-1}r^{-1}t})_r \# N^2_2(y_{t^{-1}r^{-1}r})_t
\end{eqnarray*}
for all $a, b\in A$ and $x, y\in B$, where we write $M_1=M^1_1\o M^2_1$ and $N_1=N^1_1\o N^2_1$.
\\

{\bf Proof.} (1) For all $a, b\in A$ and $x, y\in B$.
\begin{eqnarray*}
&& (a\# x)[(m\# n)_1(b\# y)]\\
&&=a_t m_1(b_{T^{-1}R^{-1}R})_{TR'}\# x_{R'} n_1(y_{T^{-1}R^{-1}T})_{Rt} \\
&&\stackrel {(2.2)}{=}(a_{tr^{-1}} m_1(b_{T^{-1}R^{-1}R})_{T})_r
 \# x_{r^{-1}r} n_1(y_{T^{-1}R^{-1}T})_{Rt} \\
&&\stackrel {(2.4)}{=}(a_{tr^{-1}T'^{-1}} m_1(b_{T^{-1}R^{-1}R}))_{Tr}
 \# x_{r^{-1}r} n_1(y_{T^{-1}R^{-1}TT'^{-1}})_{Rt} \\
&&\stackrel {(2.5)}{=}(a_{tr^{-1}T'^{-1}} m_1(b_{T^{-1}R^{-1}R}))_{rT}
 \# x_{r^{-1}r} n_1(y_{T^{-1}R^{-1}TT'^{-1}})_{tR} \\
 &&=(m_2(a_{tr^{-1}T'^{-1}}) b_{T^{-1}R^{-1}R})_{rT}
 \# x_{r^{-1}r} n_1(y_{T^{-1}R^{-1}TT'^{-1}})_{tR} \\
 &&\stackrel {(2.2)}{=}(m_2(a_{tr^{-1}T'^{-1}})_r b_{T^{-1}R^{-1}Rr'})_{T}
 \# x_{r^{-1}rr'} n_1(y_{T^{-1}R^{-1}TT'^{-1}})_{tR} \\
  &&\stackrel {(2.1)}{=}(m_2(a_{tr^{-1}T'^{-1}})_r b_{T^{-1}R^{-1}R})_{T}
 \# (x_{r^{-1}r} n_1(y_{T^{-1}R^{-1}TT'^{-1}})_{t})_R \\
 &&\stackrel {(2.3)}{=}(m_2(a_{t'^{-1}tr^{-1}T'^{-1}})_r b_{T^{-1}R^{-1}R})_{T}
 \# (x_{r^{-1}rt'^{-1}} n_1(y_{T^{-1}R^{-1}TT'^{-1}}))_{tR} \\
  &&=(m_2(a_{t'^{-1}tr^{-1}T'^{-1}})_r b_{T^{-1}R^{-1}R})_{T}
 \# (n_2(x_{r^{-1}rt'^{-1}}) y_{T^{-1}R^{-1}TT'^{-1}})_{tR} \\
 &&\stackrel {(2.5)(2.4)}{=} m_2(a_{t'^{-1}tr^{-1}T'^{-1}})_{rT} b_{T^{-1}R^{-1}RT''}
 \# (n_2(x_{r^{-1}rt'^{-1}}) y_{T^{-1}R^{-1}T''T T'^{-1}})_{Rt} \\
 &&\stackrel {(2.1)}{=} m_2(a_{t'^{-1}tr^{-1}T'^{-1}})_{rT} b_{T^{-1}R^{-1}Rr'T''}
 \# (n_2(x_{r^{-1}rt'^{-1}})_{r'} y_{T^{-1}R^{-1}T''T T'^{-1}R})_{t} \\
&&\stackrel {(2.3)}{=} m_2(a_{t'^{-1}tt''r^{-1}T'^{-1}})_{rT} b_{T^{-1}R^{-1}Rr'T''}
 \# n_2(x_{r^{-1}rt'^{-1}})_{r't} y_{T^{-1}R^{-1}T''T T'^{-1}Rt''} \\
&&\stackrel {(2.5)}{=} m_2(a_{T^{-1}R^{-1}T})_{Rt} b_r\# n_2(x_{T^{-1}R^{-1}R})_{Tr} y_t\\
 &&=(m_2(a_{T^{-1}R^{-1}T})_R\# n_2(x_{T^{-1}R^{-1}R})_T)(b\# y)\\
&& =[(m\# n)_2(a\# x)](b\# y).
\end{eqnarray*}

(2) Similarly.

This finishes the proof. $\blacksquare$
\\

{\bf Definition 2.1.10.} Let $A$ and $B$ be  multiplier Hopf algebras with two
 $k$-linear maps $R: B\o A\lr A\o B$ and $T: A\o B\lr A\o B$ as described above. One defines:
$$
 T_1^{A\# _{RT}B}=(id_{A\o B}\o R*T)\ci (\D _{A\o B}\o id_{A\o B})\ci (R^{-1})_{34}\ci (T^{-1})_{34}
 $$
 and
 $$
 T_2^{A\# _{RT}B}=(T*R\o id_{A\o B})\ci (id_{A\o B}\o \D _{A\o B})\ci (T^{-1})_{12}\ci (R^{-1})_{12}.
 $$
\\

We now prove that these candidates $T_1^{A\# _{RT}B}$ and $T_2^{A\# _{RT}B}$ defined
 by $R*T$ and $T*R$, respectively can be used to   define a good comultiplication.
\\

 {\bf Proposition 2.1.11.} Take the notations as above.
 For all $a, b, c\in A$ and $x, y, z\in B$, define
$$
\D _{RT}(a\# x)((b\# y)\o (c\# z))=T_1^{A\# _{RT}B}((a\# x)\o (b\# y))((c\# z)\o (1\# 1))
$$
and
$$
((b\# y)\o (c\# z))\D _{RT}(a\# x)=((1\# 1)\o (c\# z))T_2^{A\# _{RT}B}((a\# x)\o (b\# y)).
$$

  Then  $\D _{RT}(a\# x)$ is a two-sided multiplier of $(A\# _{RT}B)\o (A\# _{RT}B)$, where
  $\D _{RT}(a\# x)$ $=$ $\D(a)\# \D (x)$. Furthermore,  $\D _{RT}$ is coassociative
  on $A\# _{RT}B$.

 {\bf Proof.}  For all $a, b\in A$ and $x, y\in B$, we compute.
\begin{eqnarray*}
 && T_1^{A\# _{RT}B}((a\# x)\o (b\# y))\\
 &&=(id_{A\o B}\o R*T)\ci (\D _{A\o B}\o id_{A\o B})\ci (R^{-1})_{34}\ci (T^{-1})_{34}((a\# x)\o (b\# y))\\
 &&=(id_{A\o B}\o R*T)\ci (\D _{A\o B}\o id_{A\o B})((a\# x)\o (b_{t^{-1}R^{-1}}\# y_{t^{-1}R^{-1}}))\\
 &&=\sum (a_1\# x_1)\o (R*T)((a_2\# x_2)\o (b_{t^{-1}R^{-1}}\# y_{t^{-1}R^{-1}}))\\
 &&=\sum (a_1\# x_1)\o (a_2b_{t^{-1}R^{-1}R})_T\# (x_2 y_{t^{-1}R^{-1}T})_R\\
 &&\stackrel {(2.5)}{=}\sum (a_1\# x_1)\o (a_2b_{R^{-1}Rt^{-1}})_T\# (x_2 y_{t^{-1}TR^{-1}})_R\\
 &&\stackrel {(2.1)}{=}\sum (a_1\# x_1)\o (a_2b_{R^{-1}Rrt^{-1}})_T\# x_{2r} y_{t^{-1}TR^{-1}R}\\
 &&=\sum (a_1\# x_1)\o (a_2b_{rt^{-1}})_T\# x_{2r} y_{t^{-1}T}\\
 &&\stackrel {(2.4)}{=}\sum (a_1\# x_1)\o a_{2T}b_{rt^{-1}t}\# x_{2r} y_{t^{-1}tT}\\
 &&=\sum (a_1\# x_1)\o a_{2T}b_{r}\# x_{2r} y_{T}\\
 &&=\D _{RT}(a\# x)((1\# 1)\o (b\# y)).
 \end{eqnarray*}

 Similarly, by Eq.(2.2), (2.3) and (2.5) we have
$$
T_2^{A\# _{RT}B}((a\# x)\o (b\# y))=((a\# x)\o (1\# 1))\D _{RT}(b\# y).
$$

Firstly, for $a\in A, x\in B$, we prove that $\D _{RT}(a\# x)$ is a two-sided multiplier of
 $(A\# _{RT}B)\o (A\# _{RT}B)$. In fact, notice that  $\D _{RT}(a\# x)=\D _A(a)\# \D _B(x)$ and
  recall from Proposition 2.1.8  that
  for $\D _A(a)\in  M(A\o A)$ and $\D _B(x)\in  M(B \o B)$,
 we can form the multiplier $\D _A(a)\# \D _B(x)\in M((A\# _{RT}B)\o (A\# _{RT}B))$.
 We can easily check that  $\D _{RT}(a\# x)=\D _A(a)\# \D _B(x)$.
  It is clear that in the case that $A$ and $B$ are Hopf algebras,
 $\D _{RT}(a\# x)=(id_A\o \t \o id_B)(\D _A(a)\o \D _B(x))$.

 Then in order to check that $\D _{RT}$ is coassociative
  on $A\# _{RT}B$, we have to prove that
  $$
  (T_2^{A\# _{RT}B}\o id)(id\o T_1^{A\# _{RT}B})=(id\o T_1^{A\# _{RT}B})(T_2^{A\# _{RT}B}\o id).
  $$

This calculations are straightforward by making use of the expressions
$$
 T_1^{A\# _{RT}B}=(id_{A\o B}\o R*T)\ci (\D _{A\o B}\o id_{A\o B})\ci (R^{-1})_{34}\ci (T^{-1})_{34}
$$
and
$$
 T_2^{A\# _{RT}B}=(T*R\o id_{A\o B})\ci (id_{A\o B}\o \D _{A\o B})\ci (T^{-1})_{12}\ci (R^{-1})_{12}
$$
and the coassociativity of $\D _A$ and $\D _B$.

This completes the proof. $\blacksquare$
\\

We now have the main result of this section as follows.
\\

 {\bf Theorem 2.1.12.} Let $A$ and $B$ be  multiplier Hopf algebras with two
 $k$-linear maps $R: B\o A\lr A\o B$ and $T: A\o B\lr A\o B$
 such that Eq.(2.1)-(2.5) hold. If $R$ and $T$ are bijective
 such that (1): $(T*R)(x\o y\o a\o b)=0$ implies that $x\o y=0$
 and (2): $(R*T)(a\o b\o x\o y)=0$ implies that $x\o y=0$
 for all $a, x\in A$ and $b, y\in B$, and
 \begin{eqnarray*}
(2.15) && (\D _{RT}\circ R)(x\o a)=((1\o 1)\# \D _B(x))(\D _{A}(a)\# (1\o 1))\, \, \mbox{for $a\in A, x\in B$};\\
(2.16) && \D _{RT}\circ T=(T\o T)\circ \D _{RT};\\
(2.17) && (T*R)\ci (S_A\o S_B\o id_{A\o B})\circ \D _{RT} = (1\# 1)\v _{RT},\\
(2.18) && (R*T)\ci (id_{A\o B}\o S_A\o S_B)\circ \D _{RT} = (1\# 1)\v _{RT},
\end{eqnarray*}

then  $(A\# _{RT}B, \D _{RT}, \v _{RT}, S_{RT})$ is a  multiplier Hopf algebra,
 where the multiplication, $\D _{RT}$, $\v _{RT}$ and $S_{RT}$ are given as

\begin{eqnarray*}
&&(a\# b)(c\# d)=a_Tc_R\# b_Rd_T, \\
&& \D _{RT}=(id_A\o \t \o id_B)(\D _A\o \D _B),\\
&& \v _{RT}=\v _A\o \v _B,\\
&& S_{RT}= T\ci R\ci (S_B\o S_A)\t
\end{eqnarray*}
for all $a, c\in A$ and $b, d\in B$.

{\bf Proof.} We will finish the proof with the following steps:

(1) From the conditions Eq.(2.1)-(2.5) and Proposition 2.1.4, and the bijectivities of $T$ and $R$
 we have that $A\# _{RT}B$ is an associative
 algebra with non-degenerate product.

(2) From Proposition 2.1.10 we conclude that $\D _{RT}$ is coassociative
  on $A\# _{RT}B$.

(3)  We prove that $\D _{RT}: (A\# _{RT}B)\lr  M((A\# _{RT}B)\o (A\# _{RT}B))$ is a homomorphism.
  For all $a, b\in A$ and $x, y\in B$, we have
\begin{eqnarray*}
&&\D _{RT}[(a\# x)(b\# y)]=\D _{RT}(a_Tb_R\# x_Ry_T)\\
&=& \D _{A}(a_Tb_R)\# \D _{B}(x_Ry_T)\\
&=& (\D _{A}(a_T)\# (1\o 1))(\D_A(b_R)\# \D _{B}(x_R))
((1\o 1)\# \D _B(y_T))\\
&=& (\D _{A}(a_T)\# (1\o 1))(\D_{RT}\circ R) (x\o b)
((1\o 1)\# \D _B(y_T))\\
&\stackrel {(2.15)}{=}& (\D _{A}(a_T)\# (1\o 1))((1\o 1)\# \D _B(x))(\D _{A}(b)\# (1\o 1))
((1\o 1)\# \D _B(y_T))\\
&\stackrel {(2.16)}{=}& (\D _{A}(a)\# \D _B(x))(\D _{A}(b)\# \D _B(y))\\
&=&\D _{RT}(a\# x) \D _{RT}(b\# y).
\end{eqnarray*}

(4) Define the counit $\v _{RT}$ on $A\# _{RT}B$ by $\v _{RT}=\v _A\o \v _B$.
 We have to prove that for all $a, b\in A$ and $x, y\in B$

 (i) $(\v _{RT}\o id _A\o id_B)\D _{RT}(a\# x)((1\# 1)\o (b\# y))=(a\# x)(b\# y)$;

 (ii) $(id _A\o id_B\o \v _{RT})((a\# x)\o (1\# 1))\D _{RT}(b\# y)=(a\# x)(b\# y)$.

 We prove (i), the proof of (ii) is similar.
 From Proposition 2.1.10 we have that for $a, b\in A$ and $x, y\in B$
 \begin{eqnarray*}
&& \D _{RT}(a\# x)((1\# 1)\o (b\# y))\\
&&= T_1^{A\# _{RT}B}((a\# x)\o (b\# y))\\
&&=\sum (a_1\o x_1\o \sum _{ij}(a_2 b_i{}^j{}_R)_T \# (x_2 y_i{}^j{}_T)_R)\\
 \end{eqnarray*}
where  $T^{-1}(b\o y)=\sum _iy_i\o b_i$ and
  $R^{-1}(\sum _iy_i\o b_i)=\sum _{ij}b_i{}^j\o y_i{}^j$.

  So,
  \begin{eqnarray*}
&&(\v _{RT}\o id _A\o id_B) \D _{RT}(a\# x)((1\# 1)\o (b\# y))\\
&=& (\v _{RT}\o id _A\o id_B)\sum (a_1\o x_1\o \sum _{ij}(a_2 b_i{}^j{}_R)_T \# (x_2 y_i{}^j{}_T)_R)\\
&\stackrel {(2.4)}{=}& \sum _{ij}a_T b_i{}^j{}_{Rt} \# (x y_i{}^j{}_{tT})_R\\
&\stackrel {(2.1)}{=}& \sum _{ij}a_T b_i{}^j{}_{Rrt} \# x_r y_i{}^j{}_{tTR}\\
&\stackrel {(2.5)}{=}& \sum _{i}a_T b_i{}_{rt} \# x_r y_i{}_{tT}\\
&\stackrel {(2.5)}{=}& a_T b_{r} \# x_r y_{T}\\
&=& (a\# x)(b\# y).
 \end{eqnarray*}

(5) Because $T_1^{A\# _{RT}B}$ is surjective and $\D _{RT}$
 is a homomorphism, $\v _{RT}$ is a homomorphism. This can
 be proved in a similar way as in [VD1, Lemma 3.5].

(6) On $A\# _{RT}B$, define the antipode $S_{RT}= T\ci R\ci (S_B\o S_A)\t$
  which is an invertible map. We have to prove that for all $a, b\in A$
   and $x, y\in B$

(i) $m_{A\# _{RT}B}(S_{RT}\o id_A\o id_B)\D _{RT}(a\# x)((1\# 1)\o (b\# y))=\v _{RT}(a\# x)(b\# y)$;

(ii) $m_{A\# _{RT}B}( id_A\o id_B\o S_{RT})((a\# x)\o (1\# 1))\D _{RT}(b\# y)=(a\# x)\v _{RT}(b\# y)$.

 We prove (i), the proof of (ii) is similar. In fact,
 start again from the following equation:
$$
\D _{RT}(a\# x)((1\# 1)\o (b\# y))=\sum (a_1\o x_1\o \sum _{ij}(a_2 b_i{}^j{}_R)_T \# (x_2 y_i{}^j{}_T)_R)
$$
where  $T^{-1}(b\o y)=\sum _iy_i\o b_i$ and
  $R^{-1}(\sum _iy_i\o b_i)=\sum _{ij}b_i{}^j\o y_i{}^j$.

Then we have
\begin{eqnarray*}
&& m_{A\# _{RT}B}(S_{RT}\o id_A\o id_B)\D _{RT}(a\# x)((1\# 1)\o (b\# y))\\
&=& m_{A\# _{RT}B}(S_{RT}\o id_A\o id_B)\sum (a_1\o x_1\o \sum _{ij}(a_2 b_i{}^j{}_R)_T \# (x_2 y_i{}^j{}_T)_R)\\
&=&\sum \sum _{ij} S(a_1)_{R'T't}(a_2 b_i{}^j{}_R)_{Tr}\# S(b_1)_{R'T'r}
 (x_2 y_i{}^j{}_T)_{Rt}\\
&=& \sum ((S(a_1)_{T}a_{2})_R\# (S(b_1)_{R}b_{2})_T) (b\# y) \,\quad  \mbox{by (2.1)-(2.5)}\\
&\stackrel {(2.17)}{=}&\v _{RT}(a\# x)(b\# y).
 \end{eqnarray*}

(7) Because $T_1^{A\# _{RT}B}$ is surjective and $\D _{RT}$ is a homomorphism,
 $S_{RT}$ is a anti-homomorphism. The proof is similar to the proof of [VD1, Lemma 4.4].

It follows from [VD2, Proposition 2.9] that we now conclude that  $(A\# _{RT}B, \D _{RT}, $
 $\v _{RT}, S_{RT})$
 is a (regular) multiplier Hopf algebra. $\blacksquare$
\\

  {\bf Proposition 2.1.13.} Let $A$ and $B$ be multiplier
  Hopf $*$-algebras with two $k$-linear maps
  $R: B\o A\lr A\o B$ and $T: A\o B\lr A\o B$
  such that $R$ and $T$ satisfy all the conditions of Theorem 2.1.12.
  If furthermore
 \begin{eqnarray*}
(2.19) && (R(*_B\o *_A)\t)^2=id _A\o id_B,\\
(2.20) && (T(*_A\o *_B))^2=id _A\o id_B,
\end{eqnarray*}
 then   $(A\# _{RT}B, \D _{RT}, \v _{RT}, S_{RT})$ is a
 multiplier Hopf $*$-algebra.

{\bf Proof.} Straightforward. $\blacksquare$
 \\

 {\bf Proposition 2.1.14.} Let $A$ and $B$ be multiplier Hopf algebras
 as in Theorem 2.1.12. Let $\psi _A$ (resp. $\psi _B$) be a right integral
 on $A$ (resp. $B$). Then $\psi _A\o \psi _B$ is a right integral on
 the multiplier Hopf algebra $(A\# _{RT}B, \D _{RT}, \v _{RT}, S_{RT})$.

{\bf Proof.} For all $a, b\in A$ and $x, y\in B$,
\begin{eqnarray*}
 && ((\psi _A\o \psi _B\o id_{A\o B})\D _{RT}(a\# b))(x\# y)\\
 && =(\psi _A\o \psi _B\o id_{A\o B})(\D _{RT}(a\# b)((1\# 1)\o (x\# y))\\
 &&=(\psi _A\o \psi _B\o id_{A\o B})T_1^{A\# _{RT}B}((a\# x)\o (b\# y))\\
 &&=\sum (\psi _A\o \psi _B\o id_{A\o B})(a_1\# x_1)\o (a_2b_{R^{-1}Rt^{-1}})_T\# (x_2 y_{t^{-1}TR^{-1}})_R\\
 &&=\psi _A(a)\psi _B(x)(b_{R^{-1}Rt^{-1}})_T\# (y_{t^{-1}TR^{-1}})_R\\
 &&=\psi _A(a)\psi _B(x)(b\# y).
\end{eqnarray*}

This finishes the proof. $\blacksquare$
\\

{\bf Proposition 2.1.15.} Let $A$ and $B$ be multiplier Hopf algebras
 as in Theorem 2.1.12. Let $\varphi _A$ (resp. $\varphi _B$) be a left integral
 on $A$ (resp. $B$) with associated modular element $\delta _A$ (resp. $\delta _B$).
  Then the multiplier $\delta _A\# \delta _B$ is the modular
  element in $M(A\# _{RT}B)$ associated to $\varphi_A\o  \varphi _B$.
\\

{\bf Proof.} Recall from [VD2] that the modular element $\delta _A$ in $M(A)$
 for a multiplier Hopf algebra $A$  is given by $(\varphi _A\o id)\D (a)$
 when $\varphi _A(a)=1$. Now, for the
  multiplier Hopf algebra $(A\# _{RT}B, \D _{RT}, \v _{RT}, S_{RT})$, we have that modular element
  $\delta _{A\# _{RT}B}$ associated to $\varphi_A\o  \varphi _B$ is given as
  $$
  \delta _{A\# _{RT}B}=(\varphi _A\o \varphi _B\o id_{A\o B})\D _{RT}(a\# b).
  $$
We claim that $\delta _{A\# _{RT}B}=\varphi _A\o \varphi _B$

For all $a, b\in A$ and $x, y\in B$, we have
\begin{eqnarray*}
\delta _{A\# _{RT}B}(b\# y)&=& (\varphi _A\o \varphi _B\o id_{A\o B})
 (\D _{RT}(a\# b)((1\# 1)\o (x\# y))\\
  &&=(\varphi _A\o \varphi _B\o id_{A\o B})T_1^{A\# _{RT}B}((a\# x)\o (b\# y))\\
 &&=\sum (\varphi _A\o id_{A\o B})(a_1\o (a_2b_{R^{-1}Rt^{-1}})_T\# (\delta _B y_{t^{-1}TR^{-1}})_R\\
 &&=(\delta _A b_{R^{-1}Rt^{-1}})_T\# (\delta _B y_{t^{-1}TR^{-1}})_R\\
 &&\stackrel {(2.1)}{=} (\delta _Ab_{R^{-1}Rrt^{-1}})_T\# \delta _{Br} y_{t^{-1}TR^{-1}R}\\
 &&=(\delta _Ab_{rt^{-1}})_T\#\delta _{Br} y_{t^{-1}T}\\
 &&\stackrel {(2.4)}{=} \delta _{AT}b_{rt^{-1}t}\# \delta _{Br} y_{t^{-1}tT}\\
 &&= \delta _{AT}b_{r}\# \delta _{Br} y_{T}\\
 &&=(\delta _A\# \delta _B)(b\# y).
\end{eqnarray*}

This completes the proof. $\blacksquare$
\\

{\sl  2.2 Generalized smash products}
\\
\def\theequation{2. 2. \arabic{equation}}
\setcounter{equation}{0} \hskip\parindent

Let $Q$ be a regular multiplier Hopf algebra and let $A$ be a right $Q$-module algebra
 and $B$ a right $Q$-comodule algebra. Then  there are a
 $k$-linear maps $R: A\o B\lr B\o A, a\o b\mapsto \sum b_0\o (a\lhd b_{(1)})$
 for $a\in A$ and $b\in B$. By Remark 2.1.2, we may define a right smash product
  $B\# _r ^Q A$ as $k$-vector space with multiplication
\begin{equation}
(b'\# a')(b\# a)=\sum b'b_0\# (a'\lhd b_{(1)})a
\end{equation}
for $a, a'\in A, b, b'\in B$, where we use the notation $b\# a$ in place of $b\o a$
 to emphasize the new algebraic structure.
 \\

 Let $Q$ be a regular multiplier Hopf algebra and let $A$ be a left $Q$-module algebra
 and $B$ a left $Q$-comodule algebra. Then  there are another
 $k$-linear maps $R: B\o A\lr A\o B, b\o a\mapsto \sum  b_{(-1)}\rhd a \o b_0$
 for $a\in A$ and $b\in B$. Similarly, by Remark 2.1.2,
 $A\# _l ^Q B$ denotes the associative algebra structure on
 $A\o B$ given by
 \begin{equation}
(a\# b)(a'\# b')=\sum a(b_{(-1)}\rhd a')\# b_0 b'
\end{equation}
for $a, a'\in A, b, b'\in B$,
 containing $A$ and $B$ as subalgebras.
\\

We have obvious embeddings of $A$ and $B$ in the
 multiplier algebra of $B\# _r ^Q A$. And we note that
 $ B\# _r ^Q A=BA=AB$ with $ ab=\sum b_0 (a\lhd b_{(1)})$
  for $a\in A, b\in B$.
 Similar statements hold in  $A\# _l ^Q B$.
\\

{\bf Example 2.2.1} (1) Given a left coaction $\G : A\lr M(Q\o A) $ of $Q$ on $A$
 with dual right $\wh {Q}$-action $\lhd$.  One has  the right smash product
 $\wh {Q}\# _r ^{\wh Q} A$ to be the vector space $\wh {Q}\o A$ with associative
 algebra structure given for $\wh {x}, \wh {y}\in \wh {Q}$ and $a, b\in A$ by
$$
(\wh {x}\# a)(\wh {y}\# b)=\sum \wh {x} \wh {y}_{1}\# (a\lhd \wh {y}_{2})b
$$

(2) Similarly if  $\U : A\lr M(A\o Q) $ is a right coaction of $Q$ on $A$
  with dual left $\wh {Q}$-action $\rhd$,
  then the multiplication of the algebra $A\# _l ^{\wh {Q}} \wh{Q}$
   is given by
$$
(a\# \wh {x})(b\# \wh {y})=\sum a( \wh {x}_{1}\rhd b)\# \wh {x}_2 \wh {y}
$$
for $\wh {x}, \wh {y}\in \wh {Q}$ and $a, b\in A$
\\

We remark that any right  smash product $\wh {Q}\# _r ^{\wh Q} A$ can
 be identified with an associated left smash product
 $A\# _l ^{\wh {Q}^{cop}} \wh{Q}^{cop}$, where
 $\U : A\lr M(A\o Q^{op})$ is the right coaction given by
\begin{equation}
\U =(\i \o S^{-1})\ci (\t _{Q, A})\ci \G
\end{equation}
with $\t _{Q, A}: Q\o A\lr A\o Q$ being the permutation of tensor factors. In fact we have
\\

{\bf Lemma 2.2.2.}  Let $\G : A\lr M(Q\o A)$ and $\U: A\lr M(A\o Q^{op})$ be a pair of
 left and right coactions, respectively, related by (2.2.3). Then we have
 $$
 \wh {Q}\# _r ^{\wh Q} A\cong A\# _l ^{\wh {Q}^{cop}} \wh{Q}^{cop}.
 $$

{\bf Proof.} Let $\rhd : \wh {Q}^{cop}\o A\lr A$ and $\lhd : A\o \wh{Q}\lr A$
  be the left and right actions dual to $\G $ and $\U$, respectively. The equation (2.2.3)
  is equivalent to $\wh {x}\rhd a=a\lhd \wh {S^{-1}}(x)$ for all $\wh {x}\in \wh {Q}$
   and $a\in A$.

Define a linear map $f: A\# _l ^{\wh {Q}^{cop}} \wh{Q}^{cop}\lr \wh {Q}\# _r ^{\wh Q} A$ as
 $$
 f(a\# \wh {x})=\sum \wh {x}_{1}\# (a\lhd \wh {x}_{2}).
 $$
  It is easy to check that its inverse is given by
 $$
 f^{-1}(\wh {x}\# a)=\sum ( \wh {x}_{1}\rhd a)\# \wh {x}_2
 $$
 for $\wh {x}\in \wh {Q}$ and $a\in A$.

 The remaining thing is straightforward.  $\blacksquare$
\\

Let $Q, L$ be  Hopf algebras and let $A$ be a right $Q$-module algebra
 and $B$ a right $L$-comodule algebra. Let $f: L\lr Q$ be a bialgebra map.
  Then  there are a $k$-linear maps $R: A\o B\lr B\o A, a\o b\mapsto
  \sum b_0\# (a\lhd f(b_{(1)}))$
 for $a\in A$ and $b\in B$. By Remark 2.1.2, we can define a
  right smash product $B\# _r ^{(Q, L, f)} A$ as $k$-vector space with multiplication
\begin{equation}
(b'\# a')(b\# a)=\sum b'b_0\# (a'\lhd f(b_{(1)}))a
\end{equation}
for $a, a'\in A, b, b'\in B$, where we use the notation $b\# a$ in place of $b\o a$
 to emphasize the new algebraic structure.
 \\

Then  the proof of the following proposition  is straightforward.
\\

{\bf Proposition 2.2.3.} With the above notation.
  $B\# _r ^{(Q, L, f)} A$ is an associative algebra.
 \\

Similarly, let $A$ be a left $Q$-module algebra
 and $B$ a left $L$-comodule algebra.
   Let $f: L\lr Q$ be a bialgebra map. Then  there are a $k$-linear maps
    $R: B\o A\lr A\o B, b\o a\mapsto
  \sum (f(b_{(-1)})\rhd a)\# b_0$
 for $a\in A$ and $b\in B$. By Remark 2.1.2,
 $A\# _l ^{(Q, L, f)} B$ denotes the associative algebra structure on
 $A\o B$ given by
 \begin{equation}
(a\# b)(a'\# b')=\sum a (f(b_{(-1)})\rhd a')\# b_0 b'
\end{equation}
for $a, a'\in A, b, b'\in B$,
 containing $A$ and $B$ as subalgebras.
\\

 {\sl  2.3 Generalized twisted smash products}
 \\
  \def\theequation{2. 3. \arabic{equation}}
\setcounter{equation}{0} \hskip\parindent

Let $Q$ be a regular multiplier Hopf algebra. Let $A$ be a $Q$-bimodule algebra, i.e,
  $A\in {}_Q{\cal MA}_Q$  and let $B$ be a $Q$-bicomodule
 algebra, i.e.,  $B\in {}^Q{\cal MA}^Q$.
 In this section, we consider an algebra $P$ that is a generalized smash product of
  $A$ and $B$. The construction has probably been studied in [W] for Hopf algebras but not
 yet for multiplier Hopf algebras. However, the results and the arguments are very similar
 to the theory of smash products as developed in [Dr-VD-Z]. Therefore, in the following
 proposition, we do not give all the details. We concentrate on the correct statements and
 briefly indicate how things are proven.
\\

Define $B\star ^Q_rA=B\o A$ as $k$-vector space with multiplication
\begin{equation}
(b'\star a')(b\star a)=\sum b'b_0 \star (S^{-1}(b_{(-1)})\rhd a'\lhd b_{(1)}) a
\end{equation}
for $a, a'\in A, b, b'\in B$.
\\

 Similarly
 $A\star _l ^Q B$ denotes the associative algebra structure on
 $A\o B$ given by
\begin{equation}
(a\star b)(a'\star b')=\sum a(b_{(-1)}\rhd a' \lhd S^{-1}(b_{(1)}))\star b_0b'
\end{equation}
for $a, a'\in A, b, b'\in B$.
\\

{\bf Proposition 2.3.1} $B\star ^Q_rA$ and $A\star _l ^Q B$ as above are  associative
 algebras.

 {\bf Proof.} For $B\star ^Q_rA$, the twist map is given by the formula:
    $R: A\o B\lr B\o A, a\o b =\sum b_0 \o (S^{-1}(b_{(-1)})\rhd a\lhd b_{(1)})$
    for $a\in A$ and $b\in B$.   By Remark 2.1.2, $B\star ^Q_rA$ is an
    associative algebra.

    Similarly for $A\star _l ^Q B$ and the proof is finished. $\blacksquare$
\\

 The algebra $B\star ^Q_rA$ (resp. $A\star _l ^Q B$) is called a generalized right
  (resp. left) smash product, which is denoted by $P$ (resp. $U$).
 Just as in the case of smash products, we have obvious embeddings of $A$ and $B$ in the
 multiplier algebra of $P$ (resp. $U$) and if we identify these two algebras with their images in $M(P)$
 (resp. $M(U)$), we see that $P$ (resp. $U$) is the linear span of elements $ab$
 (resp. $ba$) with $a\in A$ and $b\in $ and that we
 have the commutation rules

 i) $A$ and $B$ commute,

 ii) $ab =\sum b_0 (S^{-1}(b_{(-1)})\rhd a\lhd b_{(1)})$ \\ \quad \quad \quad (resp.
 $ba=\sum (b_{(-1)}\rhd a \lhd S^{-1}(b_{(1)}))b_0$), for $a\in A, b\in B$.

Therefore we can view $P$  (resp. $U$)as the algebra generated by $A$ and $B$
 subject to these commutation rules.
\\

 {\bf Example 2.3.2.} (1) This construction reduces to well-know constructions in the following
 three special situations. If the multiplier Hopf algebra $Q$ is trivial,
 then we obtain for $P$ simply the tensor product algebra $A\o B$.
 If the left action and coaction of $Q$ on $A$ and $B$ are trivial, respectively,
 we obtain the right smash product $B\# ^Q_rA$.

 (2) Let $Q$ denote an algebraic quantum group.  Then $\wh {Q}$ is the $Q$-bimodule
  algebra. Let $B$ be a $Q$-bicomodule algebra. Then we have two generalized
   twisted smash products  $B\star ^Q_r\wh {Q}$ and $\wh {Q}\star ^Q_lB$.

(3) Let $Q$ denote a regular multiplier Hopf algebra with antipode $S^2=\i $ and
  $A$ a $Q$-bimodule algebra. Then the generalized  twisted smash product $A\star _l ^Q Q$
   is isomorphic to an ordinary smash product $A\# Q$, where the left $Q$-action
  on $Q$ is now given by $x\rightarrow a= \sum x_1\rhd a S(x_2)$,
  for all $a\in A$ and $x\in Q$.  Recall that in the original paper
  [Dr-V-Z], one developed the theory for left actions.
\\
\\
{\sl  2.4 Generalized L-R-smash products}
 \\
  \def\theequation{2. 4. \arabic{equation}}
\setcounter{equation}{0} \hskip\parindent

 The L-R-smash product was introduced and studied in a series
 of papers [P-O], with motivation  and examples coming from the
 theory of deformation quantization.
  However, we will study the case slightly different from [P-O].
\\

Let $Q$  be a multiplier Hopf algebra, $A$ a $Q$-bimodule
 algebra, i.e., $A\in {}_Q{\cal MA}_Q$  and $B$ a $Q$-bicomodule algebra,
 i.e., $B\in {}^Q{\cal MA}^Q$. Define
  $R: A\o B\lr B\o A$ and $T: B\o A\lr B\o A$
   as $R(a\o b)=\sum b_0\o a\lhd b_{(1)}$ and
  $T(b\o a)=\sum b_0\o b_{(-1)}\rhd a$, respectively,
 for all $a\in A$ and $b\in B$. Then, by Eq. (2.1.1),
 we have the generalized L-R-twisted tensor product
   $B\di ^Q_rA$  with the multiplication by
\begin{equation}
(b'\di a')(b\di a)=\sum b'_0b_0\di (a'\lhd b_{(1)})(b'_{(-1)}\rhd a)
\end{equation}
for $a, a'\in A$ and $b, b'\in B$.
\\

Similarly, Define   $R: B\o A\lr A\o B,
 b\o a\mapsto \sum (b_{(-1)}\rhd a)\o b_0$
 and $T: A\o B\lr A\o B, a\o b\mapsto \sum (a\lhd b_{(1)})\o b_0$
    for all $a\in A$ and $b\in B$. Then, by Eq. (2.1.1),
     we can define $A\di ^Q _lB=A\o B$ as $k$-vector space with multiplication
\begin{equation}
(a\di b)(c\di d)=\sum (a\lhd d_{(1)})(b_{(-1)}\rhd c)\di b_0d_0
\end{equation}
for $a, c\in A, b, d\in B$.
\\

{\bf Example 2.4.1.} (1) Let $A$ be a right $Q$-module algebra.
  Then $A$ becomes an $Q$-bimodule algebra, with left $Q$-action
  given via $\v$. In this case the multiplication of $B\di ^Q_rA$
   becomes
$$
(b'\di a')(b\di a)=\sum b'b_0\di (a'\lhd b_{(1)})a
$$
for $a, a'\in A, b, b'\in B$, hence in this case $B\di ^Q_r A$
 coincides with the generalized smash product $B\# ^Q_rA$.

(2) We note that $Q$ itself is a $Q$-bicomodule algebra.
 So, in this case, the multiplication of $Q\di ^Q_r A$
 specializes to
$$
(x\di a)(y\di b)=\sum x_2y_1\di (a\lhd y_2)(x_1\rhd b)
$$
for $a, b\in A, x, y\in Q$.
 If the left $Q$-action is trivial, then $Q\di ^Q_r A$
 coincides with the smash product $Q\# ^Q_r A$

(3) If we consider a usual Hopf algebra $Q$, and assume that $A$ is a $Q$-bimodule algebra.
  Define   $R: Q\o A\lr A\o Q, q\o a\mapsto \sum q_1\rhd a\o q_2$
 and $T: A\o Q\lr A\o Q, a\o q\mapsto \sum a\lhd q_2\o q_1$
    for all $a\in A$ and $q\in Q$. Then, by Eq. (2.1.1),
     we can define $A\di ^Q _lQ=A\o Q$ as $k$-vector space with multiplication
\begin{equation}
(a\di q)(b\di p)=\sum (a\lhd p_2)(q_1\rhd b)\di q_2p_1
\end{equation}
for $a, b\in A, q, p\in Q$. If $T(a\o q)=\sum a\lhd q_1\o q_2$
 and $R(q\o a)=\sum q_2\rhd a\o q_1$,  then the conditions (2.11) and (2.12) hold. In fact,
 we have
 \begin{eqnarray*}
&& \sum (S(a_1)_{T}a_{2})_R\# (S(q_1)_{R}q_{2})_T\\
&&= \sum ((S(a_1)\lhd (S(q_1)_{R}q_{2})_2)a_{2})_R\# (S(q_1)_{R}q_{2})_1\\
&&= \sum ((S(a_1)\lhd (S(q_1)_{R2}q_{3}))a_{2})_R\# (S(q_1)_{R1}q_{2})\\
&&= \sum S(q_1)_{1}\rhd [(S(a_1)\lhd (S(q_1)_{3}q_{3}))a_{2}]\# (S(q_1)_{2}q_{2})\\
&&= \sum [S(q_1)_{1}\rhd (S(a_1)\lhd (S(q_1)_{4}q_{3}))][S(q_1)_{2}\rhd a_{2}]
 \# (S(q_1)_{3}q_{2})\\
 &&= \sum [S(q_1)_{2}\rhd (S(a_1)\lhd (S(q_1)_{4}q_{3}))][S(q_1)_{3}\rhd a_{2}]
 \# (S(q_1)_{1}q_{2}), \\
 && \quad \mbox {by $R(q\o a)=\sum q_2\rhd a\o q_1$ }\\
&&= \sum [S(q_3)\rhd (S(a_1)\lhd (S(q_1)q_{6}))][S(q_2)\rhd a_{2}]
 \# (S(q_4)q_{5})\\
 &&= \sum [S(q_3)\rhd (S(a_1)\lhd (S(q_1)q_{4}))][S(q_2)\rhd a_{2}]
 \# 1_Q \\
 &&= \sum S(q_2)\rhd [(S(a_1)\lhd (S(q_1)q_{3}))a_{2}] \# 1_Q \\
  &&= \sum S(q_1)\rhd [(S(a_1)\lhd (S(q_2)q_{3}))a_{2}] \# 1_Q \\
  && \quad \mbox {by $T(a\o q)=\sum a\lhd q_1\o q_2$ }\\
  &&= \sum S(q)\rhd (S(a_1)a_{2}) \# 1_Q \\
&&= \v _{RT}(a\# q)(1_A\# 1_Q),
\end{eqnarray*}
and so Eq.(2.11) is obtained. Similarly, we have
\begin{eqnarray*}
\sum (a_1S(a_2)_{R})_T\# (b_1S(b_2)_T)_R= \v _{RT}(a\# b)(1_A\# 1_B).
\end{eqnarray*}

{\bf Proposition 2.4.2.} $B\di ^QA $ and $A\di ^Q_l B$  as above are associative
 algebras.

 The proof of this result is straightforward.
\\

{\bf Proposition 2.4.3.}  Let $A$ be a $Q$-bimodule algebra and let $B$ be a $Q$-bicomodule
 algebra. Then $ B\di ^Q _rA\cong B\star ^Q_r A $   and  $ A\di ^Q _lB  \cong A\star ^Q_l B$.

 {\bf Proof.} Define $\Phi : B\star ^Q _rA \lr B\di ^Q _rA$ as
  $$
  b\star a \mapsto \sum b_0 \di (b_{(-1)}\rhd a)
  $$
for all $a\in A$ and $b\in B$. First, it is not hard to check that
 the map $\Phi $ has the inverse $\Psi:  B\di ^Q _rA \lr B\star ^Q _rA$
 given by
   $$
  b\di a \mapsto \sum  b_0 \star (S^{-1}(b_{(-1)})\rhd a)
  $$
Then we check that $\Phi $ is a homomorphism as follows.
\begin{eqnarray*}
&&\Phi ((b'\star a')(b\star a))
= \sum \Phi (b'b_0 \star (S^{-1}(b_{(-1)})\rhd a'\lhd b_{(1)}) a) \\
&&\quad \quad= \sum (b'b_0)_0 \di (b'b_0)_{(-1)}\rhd [(S^{-1}(b_{(-1)})\rhd a'\lhd b_{(1)}) a] \\
&&\quad \quad=  \sum b'_0b_0 \di [(b'_{(-1)1}\rhd a'\lhd b_{(1)})][(b'_{(-1)2}b_{(-1)} \lhd a] \\
&&\quad \quad = \sum (b'_0 \di (b'_{(-1)}\rhd a'))(b_0 \di (b_{(-1)}\rhd a))=\Phi (b'\star a')\Phi (b\star a) \\
\end{eqnarray*}
for $a, a'\in A, b, b'\in B$.

Similarly for $A\di ^Q_l B $.  This  completes the proof.  $\blacksquare$
\\

{\bf Example 2.4.4.} (1) Let $Q$ be a regular multiplier Hopf algebra and let
 $B$ be a $Q$-bicomodule algebra. If $A$ is a left $Q$-module algebra regarded as an $Q$-bimodule algebra with
 trivial right $Q$-action, then $A\star ^Q_lB$ and $A\di ^Q_lB$ both coincide with $A\# ^Q_lB$,
 and the isomorphism is just the identity.

(2) If $B = Q$, the maps $\Phi : B\star ^Q _rA \lr B\di ^Q _rA$  and $\Psi:  B\di ^Q _rA \lr B\star ^Q _rA$
 become:
$$
\Phi (b\star a)=\sum b_1 \di (b_2\rhd a), \quad \mbox{and} \quad \Psi (b\di a)=\sum  b_2 \star (S^{-1}(b_1)\rhd a)
$$
for all $a\in A$ and $b\in B$.
\\

 Let now $Q$ be a multiplier Hopf algebra, $A$ a $Q$-bimodule algebra
  and $B$ a $Q$-bicomodule algebra.  Let $C$ be an algebra in
   the Yetter-Drinfeld category ${}^Q_Q{\cal YD}$ (see [De])
   that is $C$ is both a left $Q$-module algebra and a left $Q$-comodule algebra.
   These two structures are crossed via the Yetter-Drinfeld compatibility condition
    in the following sense. For all $c\in C$ and $x, y\in Q$ we require
\begin{equation}
\sum (x_1\rhd c)_{(-1)}x_2y \o (x_1\rhd c)_0=\sum x_1c_{(-1)}y\o (x_2\rhd c_0).
\end{equation}

Consider first the generalized smash product $A\# ^Q_lC$, an associative algebra.
  From Eq.(2.4.4), it follows that  $A\# ^Q_lC$ becomes a $Q$-bimodule algebra,
  with $Q$-actions
\begin{eqnarray*}
&& x\rhd (a\# c)=\sum x_1\rhd a\# x_2\rhd c\\
&& (a\# c)\lhd x= a\lhd x\o c
\end{eqnarray*}
for all $x\in Q, a\in A$ and $c\in C$, hence we may consider the algebra $(A\# ^Q_lC)\di ^Q_l B$.
 Similarly, for $B\di ^Q_r (A\# ^Q_lC)$.
\\

Meanwhile, consider the generalized smash product $C\# ^Q_lB$, an associative algebra.
 Using the condition Eq.(2.4.4), one can see that $C\# ^Q_lB$ becomes a $Q$-bicomodule algebra,
 with $Q$-coactions:
\begin{eqnarray*}
&& \G : C\# ^Q_lB\lr Q\o  (C\# ^Q_lB), \quad \G (a\# b)=\sum a_{(-1)}b_{(-1)}\o (a_0\# b_0)\\
&& \U : C\# ^Q_lB \lr (C\# ^Q_lB)\o Q, \quad \U (a\# b)=\sum (a\o b_0)\o b_{(1)}
\end{eqnarray*}
for all $a\in A$ and $b\in B$, hence we may consider the algebra $A\di ^Q_l(C\# ^Q_lB)$.
\\

{\bf Proposition 2.4.5.} We have an algebra isomorphism
$(A\# ^Q_lC)\di ^Q_l B\cong A\di ^Q_l(C\# ^Q_lB)$, given by the trivial identification.

{\bf Proof.} We compute the multiplication in $(A\# ^Q_lC)\di ^Q_l B$  as follows:
\begin{eqnarray*}
&& ((a\# c)\di b) ((a'\# c')\di b') \\
&&=\sum [((a\# c)\lhd b'_{(1)})(b_{(-1)}\rhd (a'\# c'))]\di b_0b'_0\\
&&=\sum [((a\lhd b'_{(1)} \# c))(b_{(-1)1}\rhd a'\# b_{(-1)2}\rhd c'))]\di b_0b'_0\\
&&=\sum (a\lhd b'_{(1)})(c_{(-1)}b_{(-1)1}\rhd a')\# c_0( b_{(-1)2}\rhd c')\di b_0b'_0.
\end{eqnarray*}

The multiplication in $A\di ^Q_l(C\# ^Q_lB)$ is:
\begin{eqnarray*}
&& (a\di (c\# b)) (a'\di (c'\# b')) \\
&&=\sum (a\lhd (c'\# b')_{(1)})((c\# b)_{(-1)}\rhd a')\di (c\# b)_0(c'\# b')_0\\
&&=\sum (a\lhd b'_{(1)})(c_{(-1)}b_{(-1)}\rhd a')\di (c_0\# b_0)(c'\# b'_0)\\
&&=\sum (a\lhd b'_{(1)})(c_{(-1)}b_{(-1)}\rhd a')\# c_0( b_{0(-1)}\rhd c')\di b_{00}b'_0.
\end{eqnarray*}

Hence the two multiplications coincide, completing the proof. $\blacksquare$
\\

Since the L-R-smash product coincides with the generalized smash product if
 the right $Q$-action is trivial, we also obtain:
\\

{\bf Corollary 2.4.6.} If $Q, A, B$ are as above and
 $C$  is a left $Q$-module algebra, then we have an
 algebra isomorphism $(C\# ^Q_lA)\# ^Q_lB\cong C\# ^Q_l(A\# ^Q_lB)$,
 given by the trivial identification.
\\

Let $Q$ be a regular multiplier Hopf algebra. Recall that the Drinfel'd double $D(Q)$
(generalizing the usual Drinfel'd double of a Hopf algebra) was introduced by Drabant and Van Daele
 in [Dr-VD] by a general procedure, and more explicit descriptions were
 obtained afterwards by Delvaux and Van Daele in [De-VD].
  According to one of these descriptions, the algebra structure of $D(Q)$ is
  just the twisted smash product $\wh {Q}\star Q^{cop}$. By transferring the
 whole structure of $D(Q)$ via the map $\Phi $, we can thus obtain a
 new realization of $D(Q)$, having the L-R-smash product $\wh {Q}\di Q^{cop}$
 for the algebra structure.
\\

Let $Q$ be a  multiplier Hopf algebra. From [De-VD], an invertible element $R\in M(Q\o Q)$
 is called a Drinfeld twist (or a gauge transformation) if
\begin{eqnarray}
 && (1\o R)((\i \o \D )(R))=(R\o 1)((\D \o \i)(R))\\
 && (\v \o \i )(R)=1_Q=(\i \o \v )(R).
\end{eqnarray}

If $R=R^1\o R^2\in M(Q\o Q)$ is a Drinfeld twist with inverse
 $R^{-1}=R^{-1}\o R^{-2}\in M(Q\o Q)$, then we can define a new multiplier Hopf algebra
 $Q_R$ (see, [Wa1]) with the same multiplication and counit as $Q$, for which the comultiplication
  and antipode are given by, for $x\in Q$
\begin{eqnarray}
&& \D _R(x)=R\D (x)R^{-1}, \\
&& S_R= T_RS(x)T^{-1}_R \quad \forall x\in Q
\end{eqnarray}
where $T_R=R^1S(R^2)$ is an invertible element of $M(Q)$
 with the inverse  $T_R^{-1}=S(R^{-1})R^{-2}$.
\\

{\bf Remark:} (1) It is easy to get that
$$
((\i \o \D )(R^{-1}))(1\o R^{-1})=((\D \o \i)(R^{-1}))(R^{-1}\o 1)
$$
 and
 $$
 (\v \o \i )(R^{-1})=1_Q=(\i \o \v )(R^{-1}).
$$

(2)  Let $x\in M(Q)$ be an invertible element such that $\v (x)= 1_Q$.
 If $R$ is a twist for $Q$ then so is $R^x:= R\D (x)(x^{-1}\o x^{-2})$.
  The twists $R$ and $R^x$ are said to be gauge equivalent.
\\

Let $Q$ be a multiplier Hopf algebra, $A$ a $Q$-bimodule algebra
 and $R\in M(Q\o Q)$ a Drinfeld twist. If we introduce on $A$
  another multiplication, by $a\ci b=(R^{-1}\rhd a\lhd R^1)(R^{-2}\rhd b\lhd R^2)$
 for all $a, b\in A$,  and denote this structure by ${}_{R^{-1}}A_R$, then ${}_{R^{-1}}A_R$
  is a $Q_R$-bimodule algebra, with the same  $Q$-actions as $A$. We have
\\

{\bf Proposition 2.4.6.}  Take the notations as above. Then ${}_{R^{-1}}A_R$
  is a $Q_R$-bimodule algebra.

{\bf Proof.}  We first want to show that the product in  ${}_{R^{-1}}A_R$
 is associative and non-degenerate. Furthermore, $1_Q\in A$ remains the
 unit in $M({}_{R^{-1}}A_R)$.

We compute
\begin{eqnarray*}
&&(a\ci b)\ci c=[(R^{-1}\rhd a\lhd R^1)(R^{-2}\rhd b\lhd R^2)]\ci c\\
&&=[(r^{-1}_1R^{-1}\rhd a\lhd R^1 r^1_1)(r^{-1}_2R^{-2}\rhd b\lhd R^2r^1_2)]
(r^{-2}\rhd c \lhd r^1)\\
&&=(r^{-1}\rhd a\lhd R^1)[(r^{-2}_1R^{-1}\rhd b\lhd r^1R^2_1)
(r^{-2}_2R^{-2}\rhd c \lhd r^1R^2_2)]\\
&&=(r^{-1}\rhd a\lhd R^1)[r^{-2}\rhd [(R^{-1}\rhd b\lhd r^1)(R^{-2}\rhd c \lhd r^1)]\lhd R^2] \\
&&= a\ci[(R^{-1}\rhd b\lhd r^1)(R^{-2}\rhd c \lhd r^1)]\\
&&= a\ci(b\ci c)
\end{eqnarray*}
for all $a, b, c\in A$

Then, one has to show that ${}_{R^{-1}}A_R$ is  a left $Q_R$-module algebra.

We  note that $\D _R(x)=: \sum x _{[1]}\o x_{[2]}=\sum R^1x_1R^{-1}\o R^2x_2R^{-2}$
\begin{eqnarray*}
&&x\rhd (a\ci b)=x\rhd [(R^{-1}\rhd a\lhd R^1)(R^{-2}\rhd b\lhd R^2)]\\
&&=\sum [x_1R^{-1}\rhd a\lhd R^1)][x_2R^{-2}\rhd b\lhd R^2]\\
&&=\sum (r^{-1}R^1x_1R^{-1}\rhd a\lhd r^1)(r^{-2} R^2x_2R^{-2}\rhd b\lhd r^2)\\
&&=\sum (R^1x_1R^{-1}\rhd a)\ci ( R^2x_2R^{-2}\rhd b)\\
&&=\sum (x_{[1]}\rhd a)\ci (x_{[2]}\rhd b)
\end{eqnarray*}
for all $x\in Q, a, b\in A$

Similarly, one can check that ${}_{R^{-1}}A_R$ is  a right $Q_R$-module algebra
 and that ${}_{R^{-1}}A_R$ is  a $Q_R$-bimodule. This finishes  the proof. $\blacksquare$
\\

Suppose that we have a left $Q$-comodule algebra $B$;
  then on the algebra structure of $B$ one can introduce a left $Q_R$-comodule algebra
  structure (denoted by ${}^{R^{-1}}B$ in what follows) putting the same $Q$-coaction $\G$
  as $B$. Similarly, if $C$ is a right $Q$-comodule algebra, one can
 introduce on the algebra structure of $C$ a right $Q_R$-comodule algebra structure
 (denoted by $C^{R}$ in what follows) putting  the same $Q$-coaction $\G$
  as $C$.  One may check that if $B$ is a $Q$-bicomodule algebra,
   the left and right $Q_R$-comodule algebras ${}^{R^{-1}}B$
 respectively $B^{R}$ actually
 define the structure of a $Q_R$-bicomodule algebra on $B$,
  denoted by ${}^{R^{-1}}B^R$.
\\

Then the proof of the following result is straightforward.
\\

{\bf Proposition 2.4.7.}  Take the notations as above. If $B$ is a $Q$-bicomodule algebra,
 then  ${}^{R^{-1}}B^R$ is a $Q_R$-bicomodule algebra.
\\

{\bf Proposition 2.4.8.} Let $A$ be a $Q$-bimodule algebra and
  $B$ a $Q$-bicomodule algebra. Then we have two algebra isomorphisms:
$$
B\di ^Q_r A \cong  {}^{R^{-1}}B^R \di ^Q_r {}_{R^{-1}}A_R,
  \quad \mbox {and} \quad  A\di ^Q_lB \cong {}_{R^{-1}}A_R\di ^Q_l {}^{R^{-1}}B^R
$$
given by the trivial identification.

{\bf Proof.}   We only compute the multiplication in ${}^{R^{-1}}B^R \di ^Q_r {}_{R^{-1}}A_R$.
 Similar to ${}_{R^{-1}}A_R\di ^Q_l {}^{R^{-1}}B^R$.
 \begin{eqnarray*}
&& (b'\di a')(b\di a)\\
&&=\sum b'_0b_0\di (a'\lhd b_{(1)})\ci (b'_{(-1)}\rhd a)\\
&&=\sum b'_0b_0\di (R^{-1}\rhd a'\lhd b_{(1)}R^{1}) (R^{-2}b'_{(-1)}\rhd a\lhd R^2)\\
 \end{eqnarray*}
for $a, a'\in A, b, b'\in B$.

which is the multiplication of

This finishes the proof. $\blacksquare$
\\

By Proposition 2.4.3, we have
\\

{\bf Corollary 2.4.9.} Take the notations and assumptions as above.
 We have
$$
B\star ^Q_r A \cong  {}^{R^{-1}}B^R \di ^Q_r{}_{R^{-1}}A_R,
  \quad  \mbox {and} \quad A\star ^Q_lB \cong {}_{R^{-1}}A_R\di ^Q_l{}^{R^{-1}}B^R.
$$

\section*{3. Two-times twisted tensor products}
\def\theequation{3. \arabic{equation}}
\setcounter{equation}{0} \hskip\parindent

In this section, we will study the notion of a two-times twisted tensor
 product.
\\

{\sl 3.1 Two-sided twisted tensor products}
 \\
  \def\theequation{3. 1. \arabic{equation}}
\setcounter{equation}{0} \hskip\parindent

Let $A, B, C$ be algebras  and let
 $R: C\o A\lr A\o C$ and  $T: B\o C\lr C\o B$
 be $k$-linear maps.
 Then we will write
 $R(c\o a)=a_R\o c_R=a_r\o c_r$ and $T(b\o c)=c_T\o b_T=c_t\o b_t$
 for all $a\in A, b\in B$ and $c\in C$. One defines
 two-sided twisted tensor product $A\# _RC\# _TB=A\o C\o B$
 as spaces with a new multiplication defined by the formula:
$$
m_{A\# _RC\#_TB}=(m_A\o m_C\o m_B)\circ (\i_A\o R\o T\o \i_B)
 \circ (\i_{A\o C}\o \si \o \i_{C\o B})
 $$
 or
 $$
 (a\bm c\bm b)(a'\bm c' \bm b')=aa'_R\o c_Rc'_T\o b_Tb'
$$
for all  $a, a'\in A, b, b'\in B$ and $c, c'\in C$, where
 $\si $ denotes the usual flip map on $B\o A$.
\\

Let $Q, L$ be multiplier Hopf algebras.  Let $A$ be in ${}_Q{\cal MA}$, $B\in {\cal MA}_L$, and
  let $C$ be in ${}^Q{\cal MA}^L$.
  If we define $R: C\o A\lr A\o C$ and  $T: B\o C\lr C\o B$, respectively
  by $R(c\o a)=\sum c_{(-1)}\rhd a\o c_0$
  and $T(b\o c)=c_0\o b\lhd c_{(1)}$  for all $a\in A, b\in B$ and $c\in C$,
 then we have a two-sided twisted tensor product
   $A\bm ^Q_l C\bm^L_r B$ with the multiplication given by
\begin{equation}
(a\bm c\bm b)(a'\bm c' \bm b')=\sum a(c_{(-1)}\rhd a')\bm c_0c'_0\bm (b\lhd c'_{(1)})b'
\end{equation}
for $a, a'\in A, b, b'\in B, c, c'\in C$, here
 we write $a\bm c\bm b$ for $a\o c\o b$.
 \\

 {\bf Proposition 3.1.1.} $A\bm ^Q_l C\bm^L_r B$ defined above is an associative algebra.
 \\

Note that, given $A, B$ as above, $A\o B$ becomes a $Q$-bimodule algebra,
 with $Q$-actions
$$
x \rhd (a\o b)\lhd y= x\rhd  a \o b\lhd y
$$
for all $a\in A, x, y\in Q, b\in B$.
\\

Then we have
\\

{\bf Proposition 3.1.2.} Let $A, B, C, Q$ be as above. Then
 we have the following algebra isomorphisms

(1) $ \Phi : (A\o B)\di ^L_lC\cong A\bm ^Q_l C\bm^L_r B, \quad \Phi ((a\o b)\di c)=a\bm c\bm b$;

(2) $ \Psi : C\di ^Q_r(A\o B)\cong A\bm ^Q_l C\bm^Q_r B, \quad \Psi (c\di (a\o b))=a\bm c\bm b$;

(3)  $(A\o B)\star  ^L_lC \cong (A\o B)\di ^L_lC \cong C\di ^Q_r(A\o B)\cong C\star ^Q_r(A\o B)$;

(4) $(A\# ^Q_l C)\#^L_r B= A\# ^Q_l (C\#^L_r B)= A\bm ^Q_l C\bm^L_r B$.

{\bf Proof.} (1) We do calculations as follows:
\begin{eqnarray*}
&& \Phi [((a\o b)\di c)((a'\o b')\di c')]\\
&=& \sum \Phi [((a\o b)\lhd c'_{(1)})(c_{(-1)}\rhd (a'\o b'))\di  c_0c'_0], \quad \mbox{by Eq.(2.4.3)}\\
&=&  \sum\Phi [(a\o b\lhd c'_{(1)})(c_{(-1)}\rhd a'\o b')\di  c_0c'_0]\\
&=& \sum\Phi [ a(c_{(-1)}\rhd a')\o (b\lhd c'_{(1)})b'\di  c_0c'_0]\\
&=&  \sum a(c_{(-1)}\rhd a')\bm  c_0c'_0\bm(b\lhd c'_{(1)})b'\\
&=& (a\# c\# b)(a'\bm c' \bm b')\\
&=& \Phi [((a\o b)\di c)]\Phi [((a'\o b')\di c')]
\end{eqnarray*}
for $a, a'\in A, b, b'\in B, c, c'\in C$.

(2) Similar to (1).

(3) follows (1), (2) and Proposition 2.4.3.

(4) We note that $A\# ^Q_l C$ has an induced right $Q$-comodule coaction by
 the right coaction of $Q$ on $C$, i.e., $\U: A\# ^Q_l C\lr A\# ^Q_l C \o Q$,
 $\U (a\# c)=a\# c_0\o c_{(1)}$  for all $a\in A, c\in C$. Hence,
  $(A\# ^Q_l C)\#^Q_r B= A\bm ^Q_l C\bm^Q_r B$. Similarly,
  $ A\# ^Q_l (C\#^Q_r B)= A\bm ^Q_l C\bm^Q_r B$.

This finishes the proof. $\blacksquare$
\\

{\sl 3.2 Two-sided L-R-smash products}
 \\
  \def\theequation{3. 2. \arabic{equation}}
\setcounter{equation}{0} \hskip\parindent

 Let $Q$ be a multiplier Hopf algebra, $A$ a right $Q$-comodule algebra,
  $B$ a left $Q$-comodule algebra and $C$ a $Q$-bimodule algebra.
  Define on $A\o C\o B$ a multiplication by the formula
\begin{equation}
(a\bt c\bt b)(a'\bt c' \bt b')=\sum aa'_0\bt (c\lhd a'_{(1)})(b_{(-1)}\rhd c')\bt b_0b'
\end{equation}
for $a, a'\in A, b, b'\in B, c, c'\in C$.

Then this multiplication yields an associative algebra structure, denoted by
 $A\bt ^Q_r C\bt ^Q_l B$.

 Note that, given $A, B$ as above, $A\o B$ becomes an $Q$-bicomodule algebra,
   with the following structure:
 $$
 \G : A\o B\lr M(Q\o (A\o B)), a\o b\mapsto b_{(-1)}\o (a\o b_0)=:(a\o b)_{[-1]}\o (a\o b)_{0},
 $$
 and
 $$
 \U : A\o B\lr M((A\o B)\o Q), a\o b\mapsto (a_0\o b)\o a_{(1)}=:(a\o b)_0\o (a\o b)_{[1]}.
 $$

 {\bf Proposition 3.2.1.} If $Q, A, B, C$ are as above, then we have the following
  algebra isomorphisms:

 (1) $ \phi : (A\o B)\di ^Q_r C \cong A\bt ^Q_r C\bt ^Q_l B, \quad \phi ((a\o b)\di c)=a\bt c\bt b$;

 (2) $ \psi : C\di ^Q_l (A\o B) \cong A\bt ^Q_r C\bt ^Q_l B, \quad \psi (c\di (a\o b))=a\bt c\bt b$;

 (3)  $(A\o B)\star ^Q_rC \cong (A\o B)\di ^Q_rC \cong C\di ^Q_l(A\o B) \cong C\star ^Q_l(A\o B)$;

 (4) $(A\# ^Q_r C)\# ^Q_l B= A\# ^Q_r (C\# ^Q_l B)= A\bt ^Q_r C\bt^Q_l B$.

{\bf Proof}. (1) We compute:
 \begin{eqnarray*}
&& \phi [((a'\o b')\di c')((a\o b)\di c)]\\
&&= \sum \phi [(a'\o b')_0(a\o b)_0\di (c'\lhd (a\o b)_{[1]})((a'\o b')_{[-1]}\rhd c)]\\
 &&= \sum \phi [(a'\o b'_0)(a_0\o b)\di (c'\lhd a_{(1)})(b'_{(-1)}\rhd c)]\\
 &&= \sum \phi [(a'a_0\o b'_0b)\di (c'\lhd a_{(1)})(b'_{(-1)}\rhd c)]\\
 &&= \sum a'a_0\bt (c'\lhd a_{(1)})(b'_{(-1)}\rhd c)\bt b'_0b\\
 &&=\phi ((a'\o b')\di c')\phi((a\o b)\di c).
 \end{eqnarray*}

 (2) and (3) Obvious.

 (4)  We note that $A\# ^Q_r C$ has an induced $Q$-bimodule algebra by
  the one on $C$.
 This finishes the proof. $\blacksquare$
\\

\section*{4. Twisted products  }
\def\theequation{4. \arabic{equation}}
\setcounter{equation}{0} \hskip\parindent

In this section we will introduce the notion of a Long module algebra and study the
 the twisted product of the given algebra.
\\

{\bf Definition 4.1.}  Let $Q$ be a multiplier Hopf algebra.

(1) Let $A$ be a left $Q$-module algebra
 and a left $Q$-comodule algebra. $A$ is called  a {\sl left-left $Q$-Long module algebra} if the following
 condition holds:
\begin{eqnarray}
\G (x\rhd a)(y\o 1)=\sum a_{(-1)}y\o x\rhd a_0
\end{eqnarray}
for all $x, y\in Q, a\in A$. The category of all left $Q$-Long module algebras
 is denoted by ${}^Q_Q{\cal LA}$.

 (2) Let $A$ be a left $Q$-module algebra
 and a right $Q$-comodule algebra. $A$ is called  a {\sl left-right $Q$-Long module algebra} if the following
 condition holds:
\begin{eqnarray}
\U ( x\rhd a)(1\o y)=\sum   x\rhd a_0\o a_{(1)}y
\end{eqnarray}
for all $x, y\in Q, a\in A$. The category of all left-right $Q$-Long module algebras
 is denoted by ${}_Q{\cal LA}^Q$.

(3) Let $A$ be a right $Q$-module algebra
 and a right $Q$-comodule algebra. $A$ is called  a  {\sl right-right $Q$-Long module algebra} if the following
 condition holds:
\begin{eqnarray}
\U ( a\lhd x)(1\o y)=\sum   a_0 \lhd x \o a_{(1)}y
\end{eqnarray}
for all $x, y\in Q, a\in A$. The category of all right-right $Q$-Long module algebras
 is denoted by ${}{\cal LA}^Q_Q$.

(4) Let $A$ be a right $Q$-module algebra
 and a left $Q$-comodule algebra. $A$ is called  a {\sl right-left $Q$-Long module algebra} if the following
 condition holds:
\begin{eqnarray}
\G (a\lhd x)(y\o 1)=\sum a_{(-1)}y\o a_0\lhd x
\end{eqnarray}
for all $x, y\in Q, a\in A$. The category of all right-left $Q$-Long module algebras
 is denoted by ${}^Q{\cal LA}_Q$.
\\

Let $A\in {}^Q_Q{\cal LA}$. If we define a new multiplication on $A$, by
\begin{eqnarray}
a\bu b=\sum a_0(a_{(-1)}\rhd b) \quad \forall a, b\in A
\end{eqnarray}
then this multiplication defines a new algebra structure on $A$. The product $\bu $
 is called the left twisted product.
\\

{\bf Example 4.2.}  Let $Q$ be a regular multiplier Hopf algebra with bijective antipode $S$,
 $A$  a $Q$-bimodule algebra and $B$  a $Q$-bicomodule algebra.  Then $A$
  becomes a left $Q\o Q^{op}$-module algebra with the following structure:
 $$
 (x\o y)\rhd a=x\rhd a\lhd y,
 $$
and $B$  becomes a left $Q\o Q^{op}$-comodule algebra with the following structure
$$
\G : B\lr (Q\o Q^{op})\o B, \, \, b\mapsto \sum b_{(-1)}\o S^{-1}(b_{(1)})\o b_0.
$$

It is easy to check that  $A\o B$ is
  a left $Q\o Q^{op}$-Long module algebra with the natural structure. The corresponding twisted
product on $A\o B$ is
\begin{eqnarray*}
(a\o b)\bu (a'\o b')&=&(a\o b)_0 [(a\o b)_{(-1)}\rhd (a'\o b')]\\
&=&(a\o b_0) [ b_{(-1)}\rhd a'\lhd S^{-1}(b_{(1)} \o b')]\\
&=&a (b_{(-1)}\rhd a'\lhd S^{-1}(b_{(1)})) \o b_0 b'
\end{eqnarray*}
for all $a, a'\in A$ and $b, b'\in B$, and this is exactly the multiplication of the generalized right
 twisted smash product.
\\

Similarly, for $A\in {\cal LA}^Q_Q$. If we define a new multiplication on $A$, by
\begin{eqnarray}
a\ast b=\sum (a\lhd b_{(1)})b_0 \quad \forall a, b\in A
\end{eqnarray}
then this multiplication defines a new algebra structure on $A$. The product $\ast $
 is called the right twisted product.
\\

Let $Q$ be a regular multiplier Hopf algebra with bijective antipode $S$,
 let $A$ be a $Q$-bimodule algebra and a $Q$-bicomodule algebra such that
 $A\in {}^Q_Q{\cal LA}^Q_Q$, called {\sl four sides Long module category}.  Then we define
 a new multiplication on $A$ by
 \begin{eqnarray}
a\circledast b=\sum (a_0\lhd b_{(1)})(a_{(-1)}\rhd b_0) \quad \forall a, b\in A
\end{eqnarray}
 and call it the L-R-twisted product.
 \\

Then it is easy to check that  $A$ is
  a left $Q\o Q^{op}$-Long module algebra with the natural structure
  (see Example 4.2). The corresponding twisted
 product on $A$ is
\begin{eqnarray}
a\bu b=\sum a_0(a_{(-1)}\rhd b\lhd S^{-1}(a_{(1)})) \quad \forall a, b\in A.
\end{eqnarray}

{\bf Proposition 4.3.} $(A, \circledast )$ is an associative algebra.

{\bf Proof.}  We do calculations as follows:
\begin{eqnarray*}
&&(a\circledast b)\circledast c\\
&&=\sum [(a_0\lhd b_{(1)})(a_{(-1)}\rhd b_0)]\circledast c\\
&&=\sum [(a_0\lhd b_{(1)})(a_{(-1)}\rhd b_0)]_0\lhd c_{(1)})([(a_0\lhd b_{(1)})(a_{(-1)}\rhd b_0)]_{(-1)}\rhd c_0)\\
&&=\sum [(a_0\lhd b_{(1)})_0(a_{(-1)}\rhd b_0)_0]\lhd c_{(1)})([(a_0\lhd b_{(1)})_{(-1)}(a_{(-1)}\rhd b_0)_{(-1)}]\rhd c_0)\\
&&=\sum [(a_0\lhd b_{(1)})(a_{(-1)}\rhd b_0)]\lhd c_{(1)})((a_{0(-1)}b_{0(-1)})\rhd c_0)\\
&&=\sum (a_0\lhd [b_{0(1)}c_{0(1)}](a_{(-1)}\rhd [(b_0\lhd c_{(1)})(b_{(-1)}\rhd c_0)])\\
&&=\sum (a_0\lhd [(b_0\lhd c_{(1)})_{(1)}(b_{(-1)}\rhd c_0)_{(1)}])
(a_{(-1)}\rhd [(b_0\lhd c_{(1)})_0(b_{(-1)}\rhd c_0)_0])\\
&&=\sum (a_0\lhd [(b_0\lhd c_{(1)})(b_{(-1)}\rhd c_0)]_{(1)})(a_{(-1)}\rhd [(b_0\lhd c_{(1)})(b_{(-1)}\rhd c_0)]_0)\\
&&=\sum a\circledast [(b_0\lhd c_{(1)})(b_{(-1)}\rhd c_0)]\\
&&=\sum a\circledast [(b_0\lhd c_{(1)})(b_{(-1)}\rhd c_0)]\\
&=&a\circledast (b\circledast c)
\end{eqnarray*}
for all $a, b , c\in A$.

This finishes the proof. $\blacksquare$
\\

Four sides Long module algebra category is in particular regarded as
 a left-left Long module algebra category, but in general the
 corresponding twisted products $\circledast$ and respectively $\bu $ are different.
 On the other hand, any a left-left Long module algebra category
 can be regarded as a four sides Long module algebra category
  with trivial right action and coaction, and in this case the
  corresponding twisted products coincide.
\\

{\bf Proposition 4.4.} Let $A\in {}^Q_Q{\cal LA}^Q_Q$. With
 notation as before. Then the L-R-twisted product $\circledast $
  can be obtained as a left twisting followed by a right twisting and also vice versa.

{\bf Proof.} First consider the left twisted product algebra $(A, \bu )$;
 it is easy to see that $(A, \bu )\in {\cal LA}^Q_Q$,
  and the corresponding right twisted product becomes:
\begin{eqnarray*}
a\ast b&=&\sum (a\lhd b_{(1)})\bu b_0\\
&=&\sum (a\lhd b_{(1)})_0((a\lhd b_{(1)})_{(-1)}\rhd b_0)\\
&=&\sum (a_0\lhd b_{(1)})(a_{(-1)}\rhd b_0)\\
&&=a\circledast b
\end{eqnarray*}
for all $a, b\in A$.

Similarly, one can start with the right twisted product algebra $(A, \ast)$,
 for which $(A, \ast)\in {}^Q_Q{\cal LA}$, and the corresponding
 left twisted product coincides with the L-R-twisted product.

This finishes the proof. $\blacksquare$
\\

{\bf Example 4.5.}  Let $Q$ be a regular multiplier Hopf algebra,
 $A$  a $Q$-bimodule algebra and $B$  a $Q$-bicomodule algebra.
 Take the algebra $A\o B$, which becomes a $Q$-bimodule algebra with actions
 $x\rhd (a\o b)=x\rhd a\o b$ and  $ (a\o b)\lhd x=a\lhd x\o b$,
 for all $x\in Q, a\in A, b\in B$,
 and a $Q$-bicomodule algebra, with coactions
 $\G : A\o B\lr Q\o (A\o B), \, \, a\o b\mapsto \sum b_{(-1)}\o (a\o b_0)$
 and  $\U : A\o B\lr (A\o B)\o Q, \, \, a\o b\mapsto \sum (a\o b_0) \o b_{(1)}$.
  Moreover, one checks that $A\o B\in {}^Q_Q{\cal LA}^Q_Q$,
 hence we have an L-R-twisting datum for $A\o B$.
 The corresponding L-R-twisted product is:
\begin{eqnarray*}
(a\o b)\circledast (a'\o b')&=& \sum ((a\o b)_0\lhd (a'\o b')_{(1)})((a\o b)_{(-1)}\rhd (a'\o b')_0)\\
&=& \sum ((a\o b_0)\lhd b'_{(1)})(b_{(-1)}\rhd (a'\o b'_0))\\
&=& \sum (a\lhd b'_{(1)}\o b_0)(b_{(-1)}\rhd a'\o b'_0)\\
&=& \sum (a\lhd b'_{(1)})(b_{(-1)}\rhd a')\o b_0b'_0\\
\end{eqnarray*}
for all $a, a'\in A$ and $b, b'\in B$, and this is exactly the multiplication of the generalized right
 twisted smash product.
and this is exactly the multiplication of the L-R-smash product
 $A\di ^Q _lB$.
\\

{\bf Theorem 4.6.}  With notation as above, let $Q$ be a regular multiplier  Hopf
 algebra with bijective antipode $S$. Let $A\in {}^Q_Q{\cal LA}^Q_Q$. Then
  $A$ is  a left $Q\o Q^{op}$-Long module algebra (see Example 4.2).
   Moreover, the corresponding twisted algebras $(A, \bu )$
   and $(A, \circledast )$ are isomorphic, and the isomorphism is
 defined by:
\begin{eqnarray}
&\alpha : (A, \bu ) \lr (A, \circledast ), \quad a\mapsto \sum a_0\lhd a_{(1)};\\
&\alpha ^{-1}: (A, \circledast ) \lr (A, \bu ), \quad a\mapsto \sum a_0\lhd S^{-1}(a_{(1)})
\end{eqnarray}
In particular, we obtain $ A\di ^Q _rA\cong A\star ^Q_r A $
 and  $ A\di ^Q _lA  \cong A\star ^Q_l A$.

{\bf Proof.}  We only prove that $\alpha $ is an algebra isomorphism.
 It is easy to check  that $\alpha \alpha ^{-1}=\alpha ^{-1}\alpha = id$.
 Hence we only have to check that $\alpha $ is multiplicative:
\begin{eqnarray*}
\a (a\bu b)&=&\sum \a (a_0(a_{(-1)}\rhd b\lhd S^{-1}(a_{(1)}))) \quad \mbox {by Eq.(4.8)}\\
&=& \sum (a_0(a_{(-1)}\rhd b\lhd S^{-1}(a_{(1)})))_0\lhd (a_0(a_{(-1)}\rhd b\lhd S^{-1}(a_{(1)})))_{(1)}\\
&=& \sum (a_0(a_{(-1)}\rhd b\lhd S^{-1}(a_{(1)}))_0)\lhd (a_{0(1)}(a_{(-1)}\rhd b)_{(1)})\\
&=& \sum (a_0(a_{(-1)}\rhd b_0\lhd S^{-1}(a_{(1)}))\lhd (a_{(1)} b_{(1)})
\end{eqnarray*}
and
\begin{eqnarray*}
\a (a)\circledast \a (b)&=&\sum (a_0\lhd a_{(1)})\circledast (b_0\lhd b_{(1)})\\
&=&\sum (a_0\lhd a_{(1)})\circledast (b_0\lhd b_{(1)})\\
&=&\sum  ((a_0\lhd a_{(1)})_0\lhd (b_0\lhd b_{(1)})_{(1)})
((a_0\lhd a_{(1)})_{(-1)}\rhd (b_0\lhd b_{(1)})_0)\\
&=&\sum  ((a_0\lhd a_{(1)})\lhd (b_{0(1)}))
(a_{0(-1)}\rhd (b_0\lhd b_{(1)}))\\
&=&\sum  ((a_0\lhd a_{(1)})\lhd (b_{(1)1}))
(a_{(-1)}\rhd b_0)\lhd b_{(1)2}))\\
&=&\sum  ((a_0\lhd a_{(1)})(a_{(-1)}\rhd b_0))\lhd b_{(1)}\\
\end{eqnarray*}

and the proof is finished.  $\blacksquare$
\\

\begin{center}
{\bf ACKNOWLEDGEMENT}
\end{center}

   The author would like to thank Professor A. Van Daele for his
   helpful comments on this paper. The work was partially supported
   by the NSF of   China (No. 11371088) and the NSF of China (No.11571173).

\vskip 0.8cm
\begin{center}
{\bf REFERENCES}
\end{center}

[De] L. Delvaux. Twisted tensor product of multiplier Hopf (*-)algebras.
  {\sl J. of Algebra} 269 (2003): 285-316.

 [De-VD-W] L. Delvaux, A. Van Daele \& S.H. Wang. Bicrossproducts of multiplier Hopf
  algebras. {\sl J. of Algebra} 343 (2011): 11-36.

[De-VD] Delvaux L., Van Daele A. The Drinfeld double of
 multiplier Hopf algebras. {\sl J. of Algebra} 272(1)(2004): 273-291.

[Dr-VD] Drabant B., Van Daele A.  Pairing and quantum double of
 multiplier Hopf algebras. {\sl Algebra Represent. Theory} 4(2001): 109-132.

 [Dr-VD-Z] B. Drabant, A. Van Daele \& Y. Zhang. Actions of multiplier Hopf algebra.
   {\sl Comm. Algebra} 27 (1999): 4117-4172.

[P-O] F. Panaite, F. Van Oystaeyen. L-R-smash product for (quasi-)Hopf algebras.
  {\sl J. of Algebra} 309 (2007): 168-191.

[VD1] A. Van Daele. Multiplier Hopf algebras. {\sl Trans. Am. Math. Soc.} 342(2) (1994):
  917-932.

[VD2] A. Van Daele. An algebraic framework for group duality. {\sl  Adv. in Math.} 140
  (1998): 323-366.

[VD3] A. Van Daele. Tools for working with multiplier Hopf algebras. {\sl Arabian Journal
 for Science and Engineering} 33 2C (2008): 505-527.

[VD-VK] A. Van Daele \& S. Van Keer. The Yang-Baxter and Pentagon equation. {\sl Comp.
  Math}. 91 (1994): 201-221.

[W] S. H. Wang. Doi-Koppinen Hopf bimodules are modules. {\sl Comm. in Algebra} 29(10)(2001):
 4671-4682.

[W-L]  S. H. Wang \& J. Q. Li. On twisted smash products for
 bimodule algebras and the Drinfeld double. {\sl Comm. in Algebra} 26(8)(1998):
 1435-1444.

\end{document}